\documentclass{amsart}

\usepackage{amscd}
\usepackage{epsfig}

\newtheorem{theorem}{Theorem}
\newtheorem{corollary}{Corollary}
\newtheorem{thm}{Theorem}[section]

\newtheorem{lem}[thm]{Lemma}
\newtheorem{prop}[thm]{Proposition}
\theoremstyle{definition}
\newtheorem{defn}[thm]{Definition}

\newtheorem{rem}[thm]{Remark}
\numberwithin{equation}{section}

\newcommand{\Ss}{\mathcal{S}}

\newcommand{\Ff}{\mathcal{F}}

\newcommand{\Jj}{\mathcal{J}}

\newcommand{\R}{\mathbb{R}}

\newcommand{\C}{\mathbb{C}}
\newcommand{\CP}{\mathbb{CP}}

\newcommand{\N}{{\mathbb N}}
\newcommand{\T}{{\mathbb T}}
\newcommand{\Z}{{\mathbb Z}}
\newcommand{\Q}{{\mathbb Q}}

\newcommand{\fg}{{\mathfrak{g}}}
\newcommand{\fh}{{\mathfrak{h}}}
\newcommand{\fk}{{\mathfrak{k}}}

\newcommand{\ba}{{\bf a}}
\newcommand{\bm}{{\bf m}}
\newcommand{\bo}{{\bf 1}}
\newcommand{\bba}{{[\ba]}}
\newcommand{\ha}{{\hat{a}}}

\newcommand{\QED}{\hfill {\bf Q.E.D.} \medskip}

\newcommand{\al}{{\alpha}}

\newcommand{\Om}{{\Omega}}
\newcommand{\om}{{\omega}}
\newcommand{\de}{{\delta}}

\newcommand{\ga}{{\gamma}}
\newcommand{\Ga}{{\Gamma}}

\newcommand{\la}{{\lambda}}
\newcommand{\si}{{\sigma}}
\newcommand{\La}{{\Lambda}}

\newcommand{\Si}{{\Sigma}}

\newcommand{\Diff}{{\rm Diff}}
\newcommand{\Hess}{{\rm Hess}}
\newcommand{\Det}{{\rm Det}}
\newcommand{\Ker}{{\rm Ker}}
\newcommand{\Fix}{{\rm Fix}}

\newcommand{\ov}{\overline}

\newcommand{\p}{{\partial}}

\newcommand{\tx}{\tilde{x}}
\newcommand{\tvphi}{\tilde{\varphi}}

\newcommand{\cp}{{{\CP}\,\!}}
\newcommand{\ocp}{{\ov{\CP}\,\!^2}}
\newcommand{\lp}{{{\mathbb{SP}}\,\!}}

\begin{document}

\title[K\"ahler metrics on toric orbifolds]
{K\"ahler metrics on toric orbifolds}

\author{Miguel Abreu}

\address{Fields Institute, 222 College Street, Toronto, ON M5T 3J1, Canada}
\email{mabreu@fields.utoronto.ca (until July 2001)}

\address{Departamento de Matem\'{a}tica, Instituto Superior T\'{e}cnico,
Av.Rovisco Pais, 1049-001 Lisboa, Portugal}
\email{mabreu@math.ist.utl.pt (after July 2001)}

\thanks{Partially supported by FCT (Portugal) through program POCTI and
grant POCTI/1999/MAT/33081. The author is a member of EDGE,
Research Training Network HPRN-CT-2000-00101, supported by the
European Human Potential Programme.} 



\date{\today}


\begin{abstract}
A theorem of E. Lerman and S. Tolman, generalizing a result of T. Delzant,
states that compact symplectic toric orbifolds are classified by their
moment polytopes, together with a positive integer label attached to each of
their facets. In this paper we use this result, and the existence of
``global'' action-angle coordinates, to give an effective parametrization
of all compatible toric complex structures on a compact symplectic toric
orbifold, by means of smooth functions on the corresponding moment polytope.
This is equivalent to parametrizing all toric K\"ahler metrics and
generalizes an analogous result for toric manifolds.

A simple explicit description of interesting families of extremal K\"ahler
metrics, arising from recent work of R. Bryant, is given as an application of
the approach in this paper. The fact that in dimension four these metrics
are self-dual and conformally Einstein is also discussed. This gives rise
in particular to a one parameter family of self-dual Einstein metrics
connecting the well known Eguchi-Hanson and Taub-NUT metrics.
\end{abstract}

\maketitle

\section{Introduction} \label{sec:intro}

The space of K\"ahler metrics on a K\"ahler manifold (or orbifold) can
be described in two equivalent ways, reflecting the fact that a K\"ahler
manifold is both a complex and a symplectic manifold.

From the complex point of view, one starts with a fixed complex manifold
$(M,J_0)$ and K\"ahler class $\Om\in H^{1,1}_{J_0}\cap H^2(M,\R)$,
and considers the space $\Ss (J_0,\Om)$ of all symplectic forms $\om$
on $M$ that are compatible with $J_0$ and represent the class $\Om$.
Any such form $\om\in\Ss(J_0,\Om)$ gives rise to a K\"ahler metric
$\langle\cdot,\cdot\rangle\equiv\om(\cdot,J_0\cdot)$.

The symplectic point of view arises naturally from the observation that
any two forms $\om_0,\om_1\in\Ss(J_0,\Om)$ define equivalent symplectic
structures on $M$. In fact, the family $\om_t = \om_0 + t (\om_1 - \om_0)$,
for $t\in [0,1]$, is an isotopy of symplectic forms in the same
cohomology class $\Om$, and so Moser's theorem~\cite{Mose} gives a family
of diffeomorphisms $\varphi_t : M \to M\,,\ t\in [0,1]$, such that
$\varphi_t^\ast (\om_t) = \om_0$. In particular, the K\"ahler manifold
$(M,J_0,\om_1)$ is K\"ahler isomorphic to $(M, J_1, \om_0)$, where
$J_1 = (\varphi_1)^{-1}_\ast \circ J_0 \circ (\varphi_1)_\ast$.

This means that one can also describe the space of K\"ahler metrics
starting with a fixed symplectic manifold $(M,\om_0)$ and considering
the space $\Jj (\om_0, [J_0])$ of all complex structures $J$ on $M$
that are compatible with $\om_0$ and belong to some diffeomorphism
class $[J_0]$, determined by a particular compatible complex structure
$J_0$. Any such $J\in\Jj (\om_0, [J_0])$ gives rise to a K\"ahler
metric $\langle\cdot,\cdot\rangle \equiv \om_0(\cdot,J\cdot)$.

The symplectic point of view fits into a general framework,
proposed by Donaldson in~\cite{Don1} and~\cite{Don2}, involving
the geometry of infinite dimensional groups and spaces, and
the relation between symplectic and complex quotients.
Although this framework can be useful as a guiding principle,
the symplectic point of view does not seem to be very effective
for solving specific problems in K\"ahler geometry, the reason
being that the space $\Jj (\om_0, [J_0])$ is non-linear and
difficult to parametrize. The complex point of view fairs much
better in this regard, since the space $\Ss (J_0, \Om)$ 
can be identified with an open convex subset of the linear space of smooth
functions on $M$. Indeed, the $\p\ov{\p}$-lemma asserts that given
$\om_0\in\Ss(J_0,\Om)$ any other $\om\in\Ss(J_0,\Om)$ can be
written as
\begin{equation} \label{ddbarlem}
\om = \om_0 + 2i\p\ov{\p} f\,,\ \mbox{for some}\ f\in C^\infty (M)\,.
\end{equation}
Moreover, the set of functions $f\in C^\infty(M)$ for which the
form $\om$ defined by~(\ref{ddbarlem}) is in $\Ss(J_0,\Om)$ is open
and convex.

There are however particular situations in which
the space $\Jj (\om_0, [J_0])$ admits a parametrization 
similar to the one just described for $\Ss(J_0,\Om)$, 
and the symplectic point of view can then be used very effectively.
In~\cite{Abr2} this was shown to be the case for K\"ahler toric
manifolds. In this paper we show that this can also be done for
all K\"ahler toric orbifolds, and describe an application of the
effectiveness of the symplectic approach in this context.

Let $(M,\om)$ be a symplectic toric orbifold of dimension $2n$,
equipped with an
effective Hamiltonian action $\tau:\T^n\to\Diff(M,\om)$ of the
standard (real) $n$-torus $\T^n = \R^n/2\pi\Z^n$, i.e.
$(M,\om,\tau)$ is a symplectic toric orbifold. Denote by
$\phi:M\to (\R^n)^\ast$ the moment map of such an action.
The image $P\equiv \phi(M)\subset(\R^n)^\ast$ is a convex rational
simple polytope (see Definition~\ref{def:lapo}). When $M$ is a
manifold, a theorem of Delzant~\cite{Delz} says that, up to
equivariant symplectomorphism, the polytope $P$ completely
determines the symplectic toric manifold $(M,\om,\tau)$.
In~\cite{LeTo} Lerman and Tolman generalize Delzant's theorem
to orbifolds. The result is that the polytope $P$, together
with a positive integer label attached to each of its facets,
completely determines the symplectic toric orbifold 
(see Theorem~\ref{thm:LeTo}).

The proof, in both manifold and orbifold cases, gives an
explicit construction of a canonical model for each symplectic
toric manifold, i.e. it associates to each labeled polytope $P$
an explicit symplectic toric orbifold $(M_P,\om_P,\tau_P)$ with
moment map $\phi_P:M_P\to P$ (see~\S\ref{ssec:toric}). Moreover, it
follows from the construction that $M_P$ has a canonical
$\T^n$-invariant complex structure $J_P$ compatible with $\om_P$
(see Remark~\ref{rem:kahler}).
In other words, associated to each labeled polytope 
$P\subset (\R^n)^\ast$ one has a canonical K\"ahler toric orbifold 
$(M_P,\om_P,J_P,\tau_P)$ with moment map $\phi_P:M_P\to P$.

The symplectic description of compatible toric complex structures
and K\"ahler metrics is based on the following set-up (see~\cite{Abr2}
for details). Let $\breve{P}$ denote the interior of $P$, and
consider $\breve{M}_P\subset M_P$ defined by $\breve{M}_P = \phi_P^{-1}
(\breve{P})$. One can easily check that
$\breve{M}_P$ is a smooth open dense subset of $M_P$, consisting of
all the points where the $\T^n$-action is free. It can be described as
$$
\breve{M}_P\cong \breve{P}\times \T^n =
\left\{ (x,\theta): x\in\breve{P}\subset(\R^n)^\ast\,,\ 
\theta\in\R^n/2\pi\Z^n\right\}\,,
$$
where $(x,\theta)$ are symplectic (or action-angle) coordinates for
$\om_P$, i.e.
$$
\om_P = dx\wedge d\theta = \sum_{j=1}^n dx_j \wedge d\theta_j\ .
$$

If $J$ is any $\om_P$-compatible toric complex structure on $M_P$,
the symplectic $(x,\theta)$-coordinates on $\breve{M}_P$ can be chosen
so that the matrix that represents $J$ in these coordinates has the form
\begin{equation}\label{matrixJ}
\begin{bmatrix}
\phantom{-}0\ \  & \vdots & -G^{-1} \\
\hdotsfor{3} \\
\phantom{-}G\ \  & \vdots & 0\,
\end{bmatrix}
\end{equation}
where $G=G(x)=\left[g_{jk}(x)\right]_{j,k=1}^{n,n}$ is a symmetric and
positive-definite real matrix. The integrability condition for the
complex structure $J$ is equivalent to $G$ being the Hessian of a 
smooth function $g\in C^\infty (\breve{P})$, i.e.
\begin{equation} \label{Hessg}
G = \Hess_x (g)\,,\ g_{jk}(x) = \frac{\p^2 g}{\p x_j \p x_k} (x)\,,\ 
1\leq j,k \leq n\,.
\end{equation}
Holomorphic coordinates for $J$ are given in this case by
$$
z(x,\theta) = u(x,\theta) + i v(x,\theta) = \frac{\p g}{\p x}(x)
+ i\theta\ .
$$
We will call $g$ the {\bf potential} of the compatible toric complex
structure $J$. Note that the K\"ahler metric
$\langle\cdot,\cdot\rangle = \om_P(\cdot,J\cdot)$ is given in these
$(x,\theta)$-coordinates by the matrix
\begin{equation} \label{metricG}
\begin{bmatrix}
\phantom{-}G & \vdots & 0\  \\
\hdotsfor{3} \\
\phantom{-}0 & \vdots & G^{-1}
\end{bmatrix}
\end{equation}
In particular, the induced metric on any slice of the form
$\breve{P}\times\{{\rm point}\}\subset \breve{M}_P$ is given
by the matrix $G$.

Every convex rational simple polytope $P\subset (\R^n)^\ast$ can
be described by a set of inequalities of the form
$$
\langle x, \mu_r\rangle\geq\rho_r\,,\ r=1,\ldots,d,
$$
where $d$ is the number of facets of $P$, each $\mu_r$ is a 
primitive element of the lattice $\Z^n\subset\R^n$ (the inward-pointing
normal to the $r$-th facet of P), and each $\rho_r$ is a real number.
The labels $m_r\in\N$ attached to the facets can be incorporated in
the description of $P$ by considering the affine functions
$\ell_r : (\R^n)^\ast \to \R$ defined by
$$
\ell_r (x) = \langle x, m_r\mu_r \rangle - \la_r\,\ 
\mbox{where}\ \la_r = m_r\rho_r\ \mbox{and}\ 
r=1,\ldots,d\,.
$$
Then $x$ belongs to the $r$-th facet of $P$ iff $\ell_r (x) = 0$, 
and $x\in\breve{P}$ iff $\ell_r(x) > 0$ for all $r=1,\ldots,d$.

We are now ready to state the main results of this paper.
The first is a straightforward generalization to toric orbifolds of a result
of Guillemin~\cite{Gui1}.
\begin{theorem} \label{thm1}
Let $(M_P,\om_P,\tau_P)$ be the symplectic toric orbifold associated
to a labeled polytope $P\subset (\R^n)^\ast$. Then, in suitable
symplectic $(x,\theta)$-coordinates on $\breve{M}_P\cong\breve{P}
\times \T^n$, the canonical compatible toric complex structure $J_P$ 
is of the form~(\ref{matrixJ})-(\ref{Hessg}) for a potential 
$g_P\in C^\infty (\breve{P})$ given by
$$
g_P (x) = \frac{1}{2} \sum_{r=1}^{d} \ell_r(x) \log \ell_r (x)\ .
$$
\end{theorem}
The second result provides the symplectic version of~(\ref{ddbarlem})
in this toric orbifold context, generalizing an analogous result for
toric manifolds proved in~\cite{Abr2}.
\begin{theorem} \label{thm2}
Let $J$ be any compatible toric complex structure on the symplectic
toric orbifold $(M_P,\om_P,\tau_P)$. Then, in suitable
symplectic $(x,\theta)$-coordinates on $\breve{M}_P\cong\breve{P}
\times \T^n$, $J$ is given by~(\ref{matrixJ})-(\ref{Hessg}) for a potential
$g\in C^\infty(\breve{P})$ of the form
$$
g(x) = g_P (x) + h(x)\,,
$$
where $g_P$ is given by Theorem~\ref{thm1}, $h$ is smooth on the whole
$P$, and the matrix $G=\Hess(g)$ is positive definite on $\breve{P}$
and has determinant of the form
$$
\Det(G) = \left(\de \prod_{r=1}^d \ell_r \right)^{-1}\,,
$$
with $\de$ being a smooth and strictly positive function on the whole $P$.

Conversely, any such potential $g$ determines 
by~(\ref{matrixJ})-(\ref{Hessg}) a 
complex structure on $\breve{M}_P\cong\breve{P}\times \T^n$, that
compactifies to a well-defined compatible toric complex structure $J$
on the symplectic toric orbifold $(M_P,\om_P,\tau_P)$.
\end{theorem}
Note that there is no imposed condition of $J$ being in the same
diffeomorphism class as $J_P$. The reason is that, by Theorem 9.4 
in~\cite{LeTo}, any compatible toric $J$ on $(M_P,\om_P,\tau_P)$ is
equivariantly biholomorphic to $J_P$.

Our next results describe an application of the parametrization of
compatible toric complex structures given by Theorem~\ref{thm2}.
In a recent paper~\cite{Brya} R.~Bryant studies and classifies
Bochner-K\"ahler metrics, i.e. K\"ahler metrics with vanishing
Bochner curvature. He shows in particular that these metrics always
have a very high degree of symmetry, the least symmetric ones being
of toric type. It turns out that the models for these least symmetric
Bochner-K\"ahler metrics, given by Theorem 9 in~\cite{Brya}, have
a very simple explicit description in the above symplectic framework.

For us, the most relevant geometric property of these metrics is that of 
being {\bf extremal} in the sense of Calabi (see~\S\ref{ssec:pre-ext}),
and we will construct them only as such. However, the reader should
keep in mind that these are indeed the same metrics given by
Theorem~9 in~\cite{Brya}, and hence the word ``extremal'' can be
replaced by ``Bochner-K\"ahler'' in the statements that follow.

Let $P^n_\bm$ denote the labeled simplex in $(\R^n)^\ast$ defined by the
affine functions
\begin{equation} \label{def:labl}
\ell_r(x) = m_r(1+x_r), r=1,\ldots,n\,,\ \ell_{n+1}(x) = m_{n+1} (1-\psi)\,,
\ \psi=\sum_{j=1}^n x_j\,,
\end{equation}
where $\bm = (m_1,\ldots,m_{n+1}) \in\N^{n+1}$ is a vector of positive
integer labels. The associated symplectic toric orbifold will be called
a {\bf labeled projective space} and denoted by
$(\lp^n_\bm, \om_\bm,\tau_\bm)$ ( the ``$\mathbb{S}$'' is supposed to emphasize
its $\mathbb{S}$ymplectic nature).
\begin{theorem} \label{thm3}
For any vector of labels $\bm\in\N^{n+1}$, the potential
$g\in C^\infty (\breve{P}^n_\bm)$ defined by
$$
g(x) = \frac{1}{2}\left(\sum_{r=1}^{n+1} \ell_r(x)\log\ell_r(x)
- \ell_\Si (x) \log \ell_\Si (x) \right)\,,
$$
where the $\ell_r$'s are given by~(\ref{def:labl}) and
$$
\ell_\Si (x) = \sum_{r=1}^{n+1} \ell_r (x)\,,
$$
gives rise to an extremal compatible toric complex structure on
$(\lp^n_\bm, \om_\bm,\tau_\bm)$. In other words, the metric defined
by~(\ref{metricG}) is an extremal K\"ahler metric.
\end{theorem}

As we will see in~\S\ref{ssec:weight}, there is a close relation between
labeled projective spaces $\lp^n_\bm$ and the more common weighted
projective spaces $\cp^n_\ba$. These are defined for a given vector
of positive integer weights $\ba = (a_1,\ldots,a_{n+1})\in \N^{n+1}$ as
$$
\cp^n_\ba \equiv \left(\C^{n+1}\setminus \{0\}\right)/\C^{\ast}\,,
$$
where the action of $\C^\ast = \C \setminus \{0\}$ on $\C^{n+1}$ is
given by
$$
(z_1,\ldots,z_{n+1}) \stackrel{t}{\mapsto} (t^{a_1}z_1,\ldots,t^{a_{n+1}}
z_{n+1})\,,\ t\in\C^\ast\ .
$$
The relation between $\lp^n_\bm$ and $\cp^n_\bm$ implies the following
corollary to Theorem~\ref{thm3} (see also Theorem~11 in~\cite{Brya}).
\begin{corollary} \label{cor1}
Every weighted projective space $\cp^n_\ba$ has an extremal K\"ahler
metric.
\end{corollary}

The potential $g$ of Theorem~\ref{thm3} defines an extremal K\"ahler
metric on $\breve{P}^n_\bm\times\T^n$ for any positive {\bf real} vector
of labels $\bm\in\R^{n+1}_+$. Although these do not correspond in general
to compact orbifold metrics, they do admit a natural compactification
as metrics with simple conical singularities.
\begin{theorem} \label{thm4}
Consider the smooth symplectic toric manifold
$(\lp^n_{\bf 1}\cong \cp^n, \om_{\bf 1}, \tau_{\bf 1})$ associated
to the simplex $P^n_{\bf 1}\subset (\R^n)^\ast$.
Denote by $\phi_{\bf 1}:\lp^n_{\bf 1}\to P^n_{\bf 1}$ the corresponding
moment map. Then, for any $\bm\in\R^{n+1}_+$, the extremal K\"ahler
metric~(\ref{metricG}) defined on $\breve{P}^n_{\bf 1} \times\T^n$
by the potential $g$ of Theorem~\ref{thm3}, corresponds to an extremal
K\"ahler metric on $\lp^n_{\bf 1}$ with conical singularities of
angles $2\pi/m_r$ around the pre-images $N_r\equiv \phi_{\bf 1}^{-1}(F_r)$
of each facet $F_r\subset P^n_{\bf 1}\,,\ r=1,\ldots,n+1$
(note that $N_r \cong \lp^{n-1}_{\bf 1} \cong \cp^{n-1}$).
\end{theorem}

As noted before, all extremal K\"ahler metrics of Theorem~\ref{thm4}
are in fact Bochner-K\"ahler.
In dimension four $(n=2)$ the Bochner tensor is the same as
the anti-self-dual part of the Weyl tensor, and so in this case
``Bochner-K\"ahler'' is the same as ``self-dual K\"ahler''. 
A local study and classification of these metrics in this later context
was also obtained in recent work of Apostolov and Gaudochon~\cite{ApGa}.
They show in particular that, whenever the scalar curvature $S$ is
nonzero, a self-dual K\"ahler metric is conformally Einstein with
conformal factor given by $S^{-2}$. In Section~\ref{sec:sde} we 
consider a particularly interesting family of such metrics provided
by Theorem~\ref{thm4} when $n=2$ and $\bm = (1,1,m)\,,\ m\in\R_+$.
We will see in particular that this family gives rise to a one-parameter
family of $U(2)$-invariant self-dual Einstein metrics of positive
scalar curvature, with end points the Ricci-flat Eguchi-Hanson metric
on $T\,\cp^1$ ($m=1/2$) and the also Ricci-flat Taub-NUT metric on
$\R^4$ ($m=+\infty$). We will also point out how, for a particular
discrete set of values of the parameter $m$, these metrics are
related to the ones constructed by Galicki-Lawson in~\cite{GaLa} using
quaternionic reduction.

A general discussion of the usefulness of the symplectic approach to the
construction of $U(n)$-invariant extremal K\"ahler metrics will be
given in~\cite{Abr3}.

The rest of the paper is organized as follows. 
In Section~\ref{sec:symp}, after some
necessary preliminaries on orbifolds, we give the definition and
combinatorial characterization of symplectic toric orbifolds in terms
of labeled polytopes, due to Lerman and Tolman. Labeled projective
spaces and their relation to weighted projective spaces is discussed
in~\S\ref{ssec:weight}. Theorems~\ref{thm1} and~\ref{thm2} are proved
in Section~\ref{sec:metrics}, while Theorems~\ref{thm3} and~\ref{thm4}
are proved in Section~\ref{sec:extremal}.

\section{Symplectic toric orbifolds} \label{sec:symp}

In this section, after some necessary preliminaries on orbifolds,
we give the definition and combinatorial characterization of a
symplectic toric orbifold, and discuss the family of examples given
by weighted and labeled projective spaces. 
Good references for this material are Satake~\cite{Sata} 
(for general orbifolds) and Lerman-Tolman~\cite{LeTo} 
(for symplectic orbifolds).

\subsection{Preliminaries on orbifolds} \label{ssec:prelim}

\begin{defn} \label{def:orb}
An {\bf orbifold} $M$ is a singular real manifold of dimension $n$,
whose singularities are locally isomorphic to quotient singularities
of the form $\R^n / \Ga$, where $\Ga$ is a finite subgroup of
$GL(n,\R)$ such that, for any $1\neq \ga\in\Ga$, the subspace $V_\ga
\subset\R^n$ fixed by $\ga$ has $\mbox{dim} V_\ga\leq n-2$.

For each singular point $p\in M$ there is a finite subgroup $\Ga_p
\subset GL(n,\R)$, unique up to conjugation, such that open neighborhoods
of $p$ in $M$ and $0$ in $\R^n /\Ga_p$ are homeomorphic. Such a point
$p$ is called an {\bf orbifold point} of $M$, and $\Ga_p$ the
{\bf orbifold structure group} of $p$.
\end{defn}

The condition on each nontrivial $\ga\in\Ga$ means that the singularities
of the orbifold have codimension at least two, and this makes it behave
much like a manifold. The usual definitions of vector fields, differential
forms, metrics, group actions, etc, extend naturally to orbifolds. In 
particular, a {\bf symplectic orbifold} can be defined as an orbifold $M$
equipped with a closed non-degenerate $2$-form $\om$, while a
{\bf complex orbifold} can  be defined as an orbifold $M$ equipped with
an integrable complex structure $J$. A {\bf K\"ahler orbifold} $(M,\om,J)$
is a symplectic and complex orbifold, with $\om$ and $J$ {\bf compatible}
in the sense that the bilinear form $\langle\cdot,\cdot\rangle \equiv
\om(\cdot,J\cdot)$ is symmetric and positive definite, thus defining
a {\bf K\"ahler metric} on $M$.

All orbifolds we will consider in this paper (underlying a symplectic
toric orbifold) arise through the following natural construction.
Let $Z$ be an oriented manifold and $K$ an abelian group acting
smoothly, properly and effectively on $Z$, preserving the orientation and such
that the stabilizers of points in $Z$ are always finite subgroups of $K$.
Then the quotient $M = Z/K$ is an orbifold
(the orientation preserving condition eliminates the possibility of
codimension one singularities). Its orbifold points $[p]\in M$ correspond
to points $p\in Z$ with nontrivial stabilizer $\Ga_p\subset K$,
and $\Ga_p$ is then the orbifold structure group of $[p]$.

Let $(M,\om)$ be a symplectic orbifold, and $G$ a Lie group acting
smoothly on $M$. This group action induces an infinitesimal action
of the Lie algebra $\fg$ on $M$, and for each $\xi\in\fg$
we denote by $\xi_M$ the induced vector field on $M$.
The $G$-action is said to be symplectic if it preserves $\om$,
and {\bf Hamiltonian} if it has a {\bf moment map} $\phi:M\to
\fg^\ast$, i.e. a $G$-equivariant map from $M$ to the dual
of the Lie algebra of $G$ such that
$$ 
\iota(\xi_M)\om = d\langle\xi,\phi\rangle,\ \mbox{for all}\ 
\xi\in\fg\,.
$$
When $G = S^1 = \R/2\pi\Z$, a moment map is simply given by a Hamiltonian
function $H: M\to \R^\ast \cong \fg^\ast$, whose Hamiltonian
vector field $X_H$, defined by $\iota(X_H)\om = dH$, generates
the $S^1$-action. Note that $H$ is uniquely defined up to addition by
a constant.

\subsection{Symplectic toric orbifolds} \label{ssec:toric}

\begin{defn} \label{def:torb}
A {\bf symplectic toric orbifold} is a connected $2n$-dimensional
symplectic orbifold $(M,\om)$ equipped with an effective Hamiltonian
action $\tau:\T^n \to \Diff (M,\om)$
of the standard (real) $n$-torus $\T^n = \R^n/2\pi\Z^n$.
\end{defn}

Denote by $\phi : M \to (\R^n)^\ast$ the moment map of such an action
(well-defined up to addition by a constant). When $M$ is a compact smooth
manifold, the convexity theorem of Atiyah~\cite{Atiy} and
Guillemin-Sternberg~\cite{GuS1} states that the image $P=\phi(M)
\subset (\R^n)^\ast$ of the moment map $\phi$ is the convex hull of the
image of the points in $M$ fixed by $\T^n$, i.e. a convex polytope
in $(\R^n)^\ast$. A theorem of Delzant~\cite{Delz} then says that 
the convex polytope $P\subset (\R^n)^\ast$ completely determines the 
symplectic toric manifold, up to equivariant symplectomorphisms.

In~\cite{LeTo} Lerman and Tolman generalize these two theorems to
orbifolds. While the convexity theorem generalizes word for word,
one needs more information than just the convex polytope $P$ to
generalize Delzant's classification theorem.

\begin{defn} \label{def:lapo}
A convex polytope $P$ in $(\R^n)^\ast$ is called {\bf simple} and
{\bf rational} if:
\begin{itemize}
\item[{\bf (1)}] there are $n$ edges meeting at each vertex $p$;
\item[{\bf (2)}] the edges meeting at the vertex $p$ are rational, 
i.e. each edge is of the form $p + tv_i,\ 0\leq t\leq \infty,\ 
{\rm where}\ v_i\in(\Z^n)^\ast$;
\item[{\bf (3)}] the $v_1, \ldots, v_n$ in (2) can be chosen to be a 
$\Q$-basis of the lattice $(\Z^n)^\ast$.
\end{itemize}
A {\bf facet} is a face of $P$ of codimension one. Following Lerman-Tolman,
we will say that a {\bf labeled polytope} is a rational simple convex
polytope $P\subset (\R^n)^\ast$, plus a positive integer ({\bf label}) attached
to each of its facets.

Two labeled polytopes are {\bf isomorphic} if one can be mapped to the other
by a translation, and the corresponding facets have the same integer labels.
\end{defn}

\begin{rem} In Delzant's classification theorem for compact symplectic
toric manifolds, there are no labels (or equivalently, all labels are
equal to $1$) and the polytopes that arise are
slightly more restrictive: the ``$\Q$'' in (3) is replaced by ``$\Z$''.
\end{rem}

\begin{thm}[Lerman-Tolman] \label{thm:LeTo}
Let $(M,\om,\tau)$ be a compact symplectic toric orbifold, with
moment map $\phi : M \to (\R^n)^\ast$. Then $P\equiv \phi(M)$ is a
rational simple convex polytope. For every facet $F$ of $P$, there
exists a positive integer $m_F$, the label of $F$, such that the
structure group of every $p\in \phi^{-1}(\breve{F})$ is
$\Z/m_F\Z$ (here $\breve{F}$ is the relative interior of $F$).

Two compact symplectic toric orbifolds are equivariant symplectomorphic
(with respect to a fixed torus acting on both) if and only if their
associated labeled polytopes are isomorphic. Moreover, every labeled
polytope arises from some compact symplectic toric orbifold.
\end{thm}

The proof of the last claim of this theorem is important for our
purposes. It associates to every labeled polytope $P$ a compact
symplectic toric orbifold $(M_P,\om_P,\tau_P)$, with moment map
$\phi_P : M_P \to P\subset (\R^n)^\ast$. The construction, 
generalizing Delzant's for the case of symplectic toric manifolds,
consists of a very explicit symplectic reduction.

Every labeled polytope $P\subset (\R^n)^\ast$ can be written uniquely
as
\begin{equation} \label{eq:P}
P = \bigcap_{r=1}^{d} \left\{ x\in(\R^n)^\ast\,:\ 
\ell_r(x) \stackrel{\rm def}{=} \langle x, m_r\mu_r\rangle -\la_r
\geq 0 \right\}\,,
\end{equation}
where $d$ is the number of facets, each $\mu_r$ is a primitive element
of the lattice $\Z^n \subset \R^n$ (the inward-pointing normal to the
$r$-th facet of $P$), each $m_r\in\N$ is the label attached to the
$r$-th facet of $P$, and each $\la_r$ is a real number.

Let $(e_1, \ldots, e_d)$ denote the standard basis of $\R^d$, and define
a linear map 
\begin{equation} \label{def:beta}
\beta : \R^d \to \R^n\ \ \mbox{by}\  \ \beta(e_r) = m_r \mu_r\,,\ 
r=1,\ldots,d\,. 
\end{equation}
Condition (3) of Definition~\ref{def:lapo} implies that
$\beta$ is surjective. Denoting by $\fk$ its kernel, we have short
exact sequences
$$
0 \to \fk \stackrel{\iota}{\to} \R^d \stackrel{\beta}{\to}
\R^n \to 0
\ \ \ \mbox{and its dual}\ \ \ 
0 \to (\R^n)^\ast \stackrel{\beta^\ast}{\to} (\R^d)^\ast 
\stackrel{\iota^\ast}{\to}\fk^\ast \to 0\ .
$$
Let $K$ denote the kernel of the map from $\T^d = \R^d/2\pi\Z^d$ to
$\T^n = \R^n/2\pi\Z^n$ induced by $\beta$. More precisely,
\begin{equation}\label{def:K}
K = \left\{ [\theta]\in \T^d\,:\ \sum_{r=1}^{d} \theta_r m_r \mu_r
\in 2\pi\Z^n\right\}\,.
\end{equation}
The Lie algebra of $K$ is $\fk = \Ker (\beta)$.

Consider $\R^{2d}$ with its standard symplectic form
$$
\om_0 = du\wedge dv = \sum_{r=1}^d du_r\wedge dv_r\ .
$$
We identify $\R^{2d}$ with $\C^d$ via $z_r = u_r + i v_r\,,\ 
r=1,\ldots,d$. The standard action of $\T^d$ on $\R^{2d}\cong
\C^d$ is given by
$$
\theta \cdot z = \left( e^{i\theta_1} z_1, \ldots, e^{i\theta_d} z_d\right)
$$
and has moment map
$$
\phi_{\T^d} (z_1,\ldots,z_d) = \sum_{r=1}^d \frac{|z_r|^2}{2}\, e_r^\ast 
+ \la \in (\R^d)^\ast\,,
$$
where $\la\in(\R^d)^\ast$ is an arbitrary constant. We set
$\la = \sum_r \la_r e_r^\ast$ and so
\begin{equation}\label{def:phiTd}
\phi_{\T^d} (z_1,\ldots,z_d) = \sum_{r=1}^d \left(\frac{|z_r|^2}{2} 
+ \la_r\right)e_r^\ast\in (\R^d)^\ast\ .
\end{equation}
Since $K$ is a subgroup of $\T^d$, $K$ acts on $\C^d$ with moment map
\begin{equation}\label{def:phiK}
\phi_K = \iota^\ast \circ \phi_{\T^d} =
\sum_{r=1}^d \left(\frac{|z_r|^2}{2} 
+ \la_r\right)\iota^\ast(e_r^\ast)\in \fk^\ast\ .
\end{equation}

The symplectic toric orbifold $(M_P,\om_P)$ associated to the labeled
polytope $P$ is the symplectic reduction of $\C^d$ with respect to the
$K$-action. As an orbifold it is
\begin{equation}\label{def:MPZ}
M_P = Z/K\ \ \mbox{where}\ \ Z=\phi_K^{-1}(0)\equiv\ 
\mbox{zero level set of moment map,}
\end{equation}
the symplectic structure comes from the standard one in $\R^{2d}$ (via
symplectic reduction), while the action of $\T^n \cong \T^d/K$ comes
from from the reduction of the action of $T^d$ on $Z$.

In order to verify these claims, several things need to be checked
(see \S8 of~\cite{LeTo}):
\begin{itemize}
\item[(i)] zero is a regular value of $\phi_K$ and so $Z$
is a smooth submanifold of $\R^{2d}$ of dimension $2d - (d-n) = d + n$;
\item[(ii)] with respect to the action of $K$ on $Z$, the isotropy of any
$z\in Z$ is a discrete subroup $\Ga_z$ of $K$. Hence the reduced space
$M_P = Z/K$ is a symplectic orbifold of dimension
$d+n-(d-n)=2n$;
\item[(iii)] the action of $\T^d$ on $Z$ induces an effective Hamiltonian
action of $\T^n \cong \T^d/K$ on $M_P$, whose moment map
$\phi_{\T^n}\equiv \phi_P:M_P\to (\R^n)^\ast$ has image precisely $P$;
\item[(iv)] the orbifold structure group $\Ga_{[z]}$, for any point
$[z]\in M_P$ that gets mapped by $\phi_P$ to the interior of
the $r$-th facet of $P$ (cut out by the hyperplane
$\{x\in (\R^n)^\ast :\ \ell_r (x) 
= 0\}$), is precisely $\Z/m_r\Z$.
\end{itemize}
Regarding (iii) above, recall that the moment map is apriori only
defined up to a constant. In this construction we can characterize
$\phi_P$ uniquely by requiring that it fits in the commutative
diagram
\begin{equation}  \label{def:phiP}
\begin{CD}
Z @>{\phi_{\T^d}}>> (\R^d)^\ast \\
@V{\pi}VV            @AA{\beta^\ast}A\\
M_P @>{\phi_P}>>  (\R^n)^\ast
\end{CD}
\end{equation}
where $\pi:Z\to M_P = Z/K$ is the quotient map. It is with this
normalization that $\phi_P (M_P) = P$.

\begin{rem} \label{rem:isorb}
The isotropy and orbifold structure groups of
$(M_P, \om_P, \tau_P)$ can be determined directly from the labeled
polytope $P$ (Lemma 6.6 in~\cite{LeTo}). Given $p\in M_P$, let
$\Ff(p)$ be the set of facets that contain $\phi_P(p)$, i.e.
$$
\Ff(p) = \left\{ r\in\{1,\ldots,d\} : \ell_r (\phi_P(p)) = 0\right\}\ .
$$
The isotropy group of $p$ is the subtorus $H_p\subset\T^n$ whose Lie
algebra $\fh_p$ is the linear span of the normals $\mu_r\in\R^n$,
for $r\in\Ff(p)$. The orbifold structure group $\Ga_p$ is isomorphic to
$\La_p /\hat{\La}_p$, where $\La_p\subset\fh_p$ is the lattice
of circle subgroups of $H_p$, and $\hat{\La}_p$ is the sublattice
generated by $\{m_r \mu_r\}_{r\in\Ff(p)}$.
\end{rem}

\begin{rem} \label{rem:kahler}
Note that because $(M_P, \om_P)$ is the reduction of a K\"ahler
manifold ($\C^d$ with its standard complex structure and symplectic
form) by a group action that preserves the K\"ahler structure
($K\subset U(d)$), it follows that $M_P$ comes equipped with an
invariant complex structure $J_P$ compatible with its symplectic
form $\om_P$ (see Theorem 3.5 in~\cite{GuS2}). In other words,
$(M_P,\om_P,J_P)$ is a K\"ahler toric orbifold.
\end{rem}

\subsection{Weighted and labeled projective spaces} \label{ssec:weight}

We will now discuss the family of examples of symplectic toric manifolds
given by weighted and labeled projective spaces. As we will see, these
are closely related to each other.

Consider $\C^{n+1}$ with complex coordinates $(z_1, \ldots, z_{n+1})$,
and define an action of the complex Lie group $\C^\ast = \C\setminus\{0\}$
by
\begin{equation} \label{action:c*}
(z_1,\ldots,z_{n+1}) \stackrel{t}{\mapsto}
(t^{a_1}z_1,\ldots,t^{a_{n+1}}z_{n+1})\,,\ t\in\C^\ast\,,
\end{equation}
where $a_1,\ldots,a_{n+1}$ are positive integers with highest common
divisor $1$. The {\bf weighted projective space} $\cp^n_\ba$ is defined
as the complex quotient
$$
\cp^n_\ba = \left(\C^{n+1}\setminus\{0\}\right)/\C^\ast\,,
$$
where $\ba$ denotes the vector of weights: $\ba = (a_1,\ldots,a_{n+1})$.
One checks that $\cp^n_\ba$ is a compact complex orbifold, whose orbifold
structure groups are determined in the following way. Let 
$[z]_\ba = [z_1,\ldots,z_{n+1}]_\ba$ be a point in $\cp^n_\ba$, and
let $m$ be the highest common divisor of the set of those $a_j$ for which
$z_j\neq 0$. The orbifold structure group $\Ga_{[z]_\ba}$ of $[z]_\ba$ is
isomorphic to $\Z/m\Z$. In particular, $[z]_\ba$ is a smooth point of
$\cp^n_\ba$ if and only if $m=1$. Since we assumed that the highest
common divisor of all the $a_j$'s  is $1$, this means that any point
$[z]_\ba = [z_1,\ldots,z_{n+1}]_\ba \in \cp^n_\ba$, with all $z_j\neq 0$,
is a smooth point. Note also that $\cp^n_{\bf 1}$ is the usual complex
projective space $\cp^n$, and we will omit the subscript ${\bf 1}$ when
referring to it.

There is a natural holomorphic map $\pi_\ba:\cp^n_\ba\to\cp^n$ defined by
$$
\pi_\ba \left( [z_1,\ldots,z_{n+1}]_\ba \right)
\mapsto [z_1^{\ha_1},\ldots,z_{n+1}^{\ha_{n+1}}]\,,
$$
where $\ha_j$ denotes the product of all the weights except the $j$-th
one:
$$
\ha_r = \prod_{k=1,k\neq r}^{n+1} a_k\ .
$$
The map $\pi_\ba$ factors through the quotient of $\cp^n_\ba$ by the
following finite group action. Let $\ha = \prod_{k=1}^{n+1} a_k$
and consider the finite group $\Ga_\ba$ defined by
$$
\Ga_\ba = \left(\Z_{\ha_1}\times\cdots\times\Z_{\ha_{n+1}}\right)/
\Z_\ha\,,
$$
where 
\begin{eqnarray*}
\Z_\ha & \hookrightarrow & \Z_{\ha_1}\times\cdots\times\Z_{\ha_{n+1}}\\
\zeta & \mapsto & \left(\zeta^{a_1}\,,\  \ldots\ ,\,\zeta^{a_{n+1}}\right)
\end{eqnarray*}
(here $\Z_q\equiv \Z/q\Z$ is identified with the group of $q$-th roots
of unity in $\C$).
$\Ga_\ba$ acts on $\cp^n_\ba$ via
$$
[\eta]\cdot [z]_\ba = [\eta_1 z_1,\ldots,\eta_{n+1} z_{n+1}]_\ba\,,\ 
\mbox{for all}\ [\eta]\in\Ga_\ba\,,\ [z]_\ba\in\cp^n_\ba\,,
$$
and one checks easily that
$$
\pi_\ba ([z]_\ba) = \pi_\ba ([z']_\ba)\ \mbox{iff}\ 
[z']_\ba = [\eta]\cdot [z]_\ba\ \mbox{for some}\ [\eta]\in\Ga_\ba\,.
$$
Hence we have the following commutative diagram:
\begin{eqnarray}\label{def:[pia]}
\cp^n_\ba & \stackrel{\pi_\ba}{\longrightarrow} & \cp^n \nonumber\\
\searrow & & \nearrow_{[\pi_\ba]}\\
 & \cp^n_\ba /\Ga_\ba & \nonumber
\end{eqnarray}
The action of $\Ga_\ba$ is free on $\breve{\cp}^n_\ba
= \left\{ [z_1,\ldots,z_{n+1}]_\ba\in\cp^n_\ba : z_j\neq 0\ 
\mbox{for all}\ j\right\}$. In particular, $\pi_\ba$ has degree
$|\Ga_\ba | = \left(\ha\right)^{n-1}$. It is also clear that, if
$\ba$ is a nontrivial weight vector, the $\Ga_\ba$-action has nontrivial
isotropy at some points in $\cp^n_\ba \setminus \breve{\cp}^n_\ba$,
and so $\cp^n_\ba / \Ga_\ba$ has a nontrivial orbifold structure.
The bijection $[\pi_\ba]$, although a biholomorphism between 
$\cp^n_\ba / \Ga_\ba$ and the standard $\cp^n$, is obviously not an
orbifold isomorphism between $\cp^n_\ba / \Ga_\ba$ and the smooth
$\cp^n$. We will look at $[\pi_\ba]$ as inducing an orbifold structure
on $\cp^n$ isomorphic to $\cp^n_\ba / \Ga_\ba$.
\begin{defn} \label{def:orbproj}
The {\bf orbifold projective space} $\cp^n_{[\ba]}$ is defined as the
finite quotient
$$
\cp^n_{[\ba]} \stackrel{\rm def}{=} \cp^n_\ba / \Ga_\ba
\stackrel{[\pi_\ba]}{=}\mbox{``orbifold''}\ \cp^n\ .
$$
\end{defn}
\begin{rem} \label{rem:cover}
Once the orbifold structures are taken into account, the map
$\pi_\ba:\cp^n_\ba\to\cp^n_{[\ba]}$ is an orbifold covering map of
degree $(\ha)^{n-1}$. In particular, any orbifold geometric structure
(symplectic, K\"ahler, etc) on $\cp^n_{[\ba]}$ lifts through $\pi_\ba$
to an orbifold geometric structure on $\cp^n_\ba$. For our purposes
it is then enough, and, as we will see, also more convenient, to
work with $\cp^n_{[\ba]}$.
\end{rem}

In order to better understand $\cp^n_\bba$, in particular its orbifold
structure groups and symplectic description in terms of labeled
polytopes, it is useful to go back to $\C^{n+1}$ and consider a finite
extension of the $\C^\ast$-action defined by~(\ref{action:c*}).

Let $K^\C_\ba$ be the complex Lie group defined by 
\begin{equation} \label{def:Ka}
K^\C_\ba = \left( \Z_{\ha_1}\times\cdots\times\Z_{\ha_{n+1}}\times
\C^\ast\right)/\Z_\ha
\end{equation}
where 
\begin{eqnarray*}
\Z_\ha & \hookrightarrow & \Z_{\ha_1}\times\cdots\times\Z_{\ha_{n+1}}
\times\C^\ast \\
\zeta & \mapsto & \left(\zeta^{a_1}\,,\  \ldots\ ,\,\zeta^{a_{n+1}}
\,,\ \zeta^{-1}\right)\ .
\end{eqnarray*}
$K^\C_\ba$ acts effectively on $\C^{n+1}$ via
\begin{equation}\label{action:K}
[({\bf \eta},t)]\cdot z = (\eta_1 t^{a_1} z_1, \ldots ,
\eta_{n+1} t^{a_{n+1}} z_{n+1})\,,\ \mbox{for all}\ [({\bf \eta},t)]\in K^\C_\ba\,,
\ z\in\C^{n+1}\,.
\end{equation}
Because of the exact sequence
$$
\begin{array}{ccccccccc}
1 & \to & \C^\ast & \hookrightarrow & K^\C_\ba & \to & \Ga_\ba & \to & 1 \\
  &     & t & \mapsto & [({\bf 1},t)] &           &         &     &   \\
  &     &   &         & [({\bf \eta},t)] & \mapsto   & [{\bf \eta}] & &  
\end{array}
$$
we have that
\begin{equation} \label{eq:orbproj}
\left( \C^{n+1}\setminus\{ 0\}\right)/ K^\C_\ba \cong
\left[\left( \C^{n+1}\setminus\{ 0\}\right)/\C^\ast\right]/\Ga_\ba =
\cp^n_\ba / \Ga_\ba = \cp^n_{[\ba]}\ .
\end{equation}
Hence, the orbifold structure of $\cp^n_{[\ba]}$ can be described
directly from the different isotropy subgroups of the
$K^\C_\ba$-action on $\C^{n+1}$ (see Lemma~\ref{lem:iso} below).

We will now give the symplectic description, in terms of labeled polytopes,
for the orbifold projective spaces $\cp^n_{[\ba]}$. Recall that the
polytope corresponding to $\cp^n$, with symplectic (K\"ahler) form
in the same cohomology class as the first Chern class,
is the simplex $P_{\bf 1}^n$ in $(\R^n)^\ast$ defined by
\begin{equation} \label{eq:stsimp}
P_{\bf 1}^n = \bigcap_{r=1}^{n+1} \left\{ x \in (\R^n)^\ast :
\ell_r(x)\equiv \langle x,m_r\mu_r\rangle - \la_r \geq 0 \right\}\,,
\end{equation}
where: $\la_r=m_r=1\,,\ r=1,\ldots,n+1$; $\mu_r = e_r\,,\ r=1,\ldots,n$, and
$\mu_{n+1} = - \sum_{j=1}^n e_j$. Here $(e_1,\ldots,e_n)$ denotes the 
standard basis of $\R^n$.
\begin{figure}[h]
\center{\epsfig{file=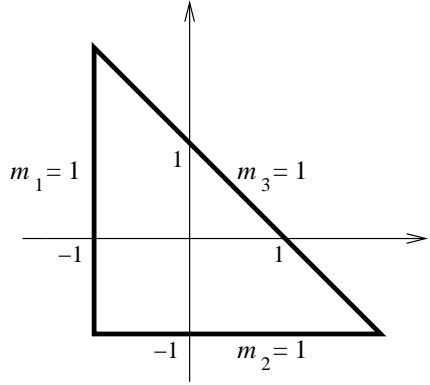}}
\caption{The polytope $P_{\bf 1}^2$ corresponding to $\cp^2$.}
\end{figure}

From \S\ref{ssec:toric} we know that a facet of a labeled polytope has
label $m\in\N$ if and only if the orbifold structure group of the
points that are mapped to its relative interior, via the moment map,
is $\Z_m \equiv \Z/m\Z$. In the case of $\cp^n$, the pre-image of the
$r$-th facet
$$
F_r = \left\{ x\in P_{\bf 1}^n : \ell_r (x) = 0 \right\}
$$
is
$$
N_r = \left\{ [z_1,\ldots,z_{n+1}]\in\cp^n : z_r = 0 \right\}\,,
$$
while the pre-image of its interior $\breve{F}_r$ is
$$
\breve{N}_r = \left\{ [z_1,\ldots,z_{n+1}]\in N_r :
z_k \neq 0\ \mbox{for all}\ k\neq r\right\}\ .
$$
In $\cp^n_{[\ba]}$ this corresponds to
$$
\breve{N}_{[\ba],r} = \left\{ [z_1,\ldots,z_{n+1}]_{[\ba]} \in
\cp^n_{[\ba]} : z_r = 0 \ \mbox{and}\ z_k\neq 0\ \mbox{for all}\ 
k\neq r \right\}\,.
$$
\begin{lem} \label{lem:iso}
The orbifold structure group $\Ga_{[z]_{[\ba]}}$ of any point
$[z]_{[\ba]}\in \breve{N}_{[\ba],r}\subset\cp^n_{[\ba]}$
is isomorphic to $\Z_{m_r}$ where
$$
m_r = \ha_r = \prod_{k=1,k\neq r}^{n+1} a_k\ .
$$ 
\end{lem}
\proof{}
Because of (\ref{eq:orbproj}), the orbifold structure group 
$\Ga_{[z]_{[\ba]}}$ of any point $[z]_{[\ba]}\in 
\breve{N}_{[\ba],r}$ is the isotropy of the $K^\C_\ba$-action
at any point $z=(z_1,\ldots,z_{n+1})\in\C^{n+1}$ with $z_r = 0$, and
$z_k\neq 0$ for all $k\neq r$. It follows from~(\ref{action:K}) that such
an isotropy subgroup is given by the elements $[({\bf \eta},t)]\in
K^\C_\ba$ such that $\eta_k = t^{-a_k}$, for all $k\neq r$.
Since $\eta_k \in \Z_{\ha_k}$, this implies that $t\in\Z_\ha \subset
\C^\ast$, and so
$$
\Ga_{[z]_{[\ba]}} \cong \left(\Z_{\ha_r}\times \Z_\ha\right) /
\left( (\zeta^{a_r},\zeta^{-1})\,,\ \zeta\in\Z_\ha\right)\ .
$$
The right-hand side is isomorphic to $\Z_{\ha_r}$ via the map
$$
[(\eta_r,\zeta)]\mapsto \eta_r \cdot \zeta^{a_r}\,,\ 
\eta_r \in \Z_{\ha_r}\,,\ \zeta\in\Z_\ha\ .
$$
\QED

The natural candidate for labeled polytope corresponding to
$\cp^n_{[\ba]}$ is then the labeled simplex $P_{[\ba]}^n$
in $(\R^n)^\ast$ defined by
\begin{equation} \label{eq:labsimp}
P_{[\ba]}^n = \bigcap_{r=1}^{n+1} \left\{ x \in (\R^n)^\ast :
\ell_r(x)\equiv \langle x,m_r\mu_r\rangle - \la_r \geq 0 \right\}\,,
\end{equation}
where $\la_r=m_r=\ha_r\,,\ r=1,\ldots,n+1$, 
and the $\mu_r$'s are as in~(\ref{eq:stsimp}).
\begin{figure}[h]
\center{\epsfig{file=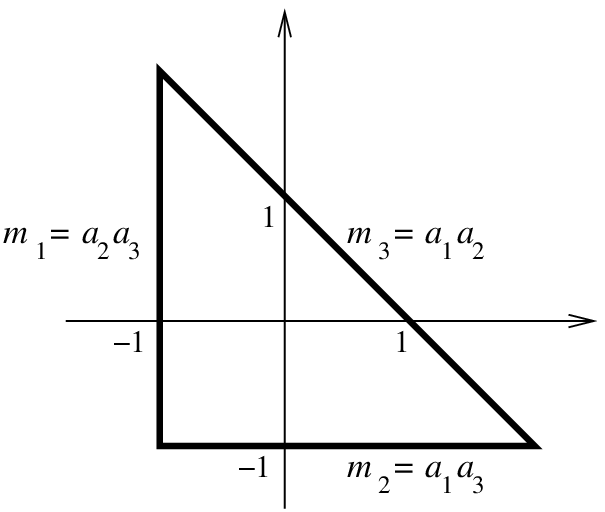}}
\caption{The labeled simplex $P_{\bba}^2$ corresponding to $\cp^2_\bba$.}
\end{figure}
\begin{prop} \label{prop:orbproj}
The compact K\"ahler toric orbifold $(M_{[\ba]}, \om_{[\ba]}, J_{[\ba]})$,
associated to the labeled polytope $P_{[\ba]}^n$ via the construction of
\S\ref{ssec:toric}, is isomorphic as a complex toric orbifold to
$\cp^n_{[\ba]}$.
\end{prop}
\proof{}
With respect to the standard basis of $\R^{n+1}$ and $\R^n$, the linear map
$\beta:\R^{n+1}\to\R^n$ defined by~(\ref{def:beta}) is given by the matrix
$$
\begin{bmatrix}
m_1 & 0 & \dots & 0 & -m_{n+1} \\
0 & \ddots & \ddots & \vdots & \vdots \\
\vdots & \ddots & \ddots & 0 & \vdots\\
0 & \dots & 0 & m_n & -m_{n+1}
\end{bmatrix}
=
\begin{bmatrix}
\ha_1 & 0 & \dots & 0 & -\ha_{n+1} \\
0 & \ddots & \ddots & \vdots & \vdots \\
\vdots & \ddots & \ddots & 0 & \vdots\\
0 & \dots & 0 & \ha_n & -\ha_{n+1}
\end{bmatrix}
$$
Using multiplicative notation, the kernel $K_\ba \subset \T^{n+1}$ of the
induced map $\beta:\T^{n+1}\to\T^n$ is then given by
$$
K_\ba = \left\{ (e^{i\theta_1},\ldots,e^{i\theta_{n+1}})\in\T^{n+1}\,:\ 
e^{i\ha_1\theta_1} = \cdots = e^{i\ha_{n+1}\theta_{n+1}}\,\right\}\ .
$$
$K_\ba$ acts on $\C^{n+1}$ as a subgroup of $\T^{n+1}$, and 
from~(\ref{def:MPZ}) and Remark~\ref{rem:kahler} we have that
$(M_\bba,\om_\bba,J_\bba)$ is the K\"ahler reduction
\begin{equation} \label{sympred}
M_\bba = \phi^{-1}_{K_\ba} (0)/ K_\ba\,,
\end{equation}
where $\phi^{-1}_{K_\ba}$ is the moment map defined by~(\ref{def:phiK}).

One easily checks that $K_\ba$ is isomorphic to the Lie group
$$
\left( \Z_{\ha_1}\times\cdots\times\Z_{\ha_{n+1}}\times
\T^1\right)/\Z_\ha
$$
where 
\begin{eqnarray*}
\Z_\ha & \hookrightarrow & \Z_{\ha_1}\times\cdots\times\Z_{\ha_{n+1}}
\times\T^1 \\
\zeta & \mapsto & \left(\zeta^{a_1}\,,\  \ldots\ ,\,\zeta^{a_{n+1}}
\,,\ \zeta^{-1}\right)\,,
\end{eqnarray*}
the isomorphism being given explicitly by
$$
\left[(\eta, e^{i\theta})\right] \mapsto
(\eta_1 e^{i a_1\theta}, \ldots ,\eta_{n+1} e^{i a_{n+1} \theta})
\in K_\ba \subset \T^{n+1}\,. 
$$
This means that the complex Lie group $K_\ba^\C$ defined by~(\ref{def:Ka})
is the complexification of $K_\ba$, and by~(\ref{eq:orbproj}) we know that
\begin{equation} \label{compquot}
\cp^n_\bba \cong (\C^{n+1}\setminus \{ 0\})/K_\ba^\C\ .
\end{equation}

The statement of the proposition now follows from a general principle
that gives an identification between the K\"ahler reduction~(\ref{sympred})
and the complex quotient~(\ref{compquot}). A good reference in our context
is the appendix to~\cite{Gui2}.
\QED

The construction of \S\ref{ssec:toric} applies of course to any
labeled polytope, and hence to any labeled simplex
\begin{equation} \label{def:labsimp}
P_{\bm}^n = \bigcap_{r=1}^{n+1} \left\{ x \in (\R^n)^\ast :
\ell_r(x)\equiv \langle x,m_r\mu_r\rangle - \la_r \geq 0 \right\}\,,
\end{equation}
with arbitrary $m_r=\la_r\in\N\,,\ r=1,\ldots,n+1$, and the
$\mu_r$'s again as in~(\ref{eq:stsimp}).
\begin{defn} \label{def:labproj}
Given an arbitrary vector $\bm = (m_1,\ldots,m_{n+1})$ of positive
integer labels, we define the {\bf labeled projective space}
$(\lp^n_\bm,\om_\bm,\tau_\bm)$ as the symplectic toric orbifold
associated to the labeled simplex $P_\bm^n \subset (\R^n)^\ast$ defined
by~(\ref{def:labsimp}).
\end{defn}
\begin{rem} \label{rem:spcp}
It follows from Remark~\ref{rem:kahler} that any labeled projective
space $(\lp^n_\bm,\om_\bm,\tau_\bm)$ comes equipped with a
``canonical'' compatible toric complex structure $J_\bm$.
Theorem~9.4 in~\cite{LeTo} implies that as a complex toric variety,
not only with respect to $J_\bm$ but also with respect to any
toric complex structure $J$ compatible with $\om_\bm$, 
$\lp^n_\bm$ is equivariantly biholomorphic to $\cp^n$. The
biholomorphism $[\pi_\ba]:\cp^n_\bba\to\cp^n$ defined by~(\ref{def:[pia]})
is just a particular explicit instance of this more general fact.
\end{rem}
\begin{rem} \label{rem:norm}
In Definition~\ref{def:labproj} we have normalized all labeled simplices
$P^n_\bm$ by the conditions $m_r = \la_r\,,\ r=1,\ldots,n+1$, which
amounts to the fact that the underlying simplex is always the same
$P^n_\bo\subset(\R^n)^\ast$. This also means that the cohomology class
of $\om_\bm$ in $H^2(\lp^n_\bm)$ is apriori fixed. One can allow for
an arbitrary positive scaling of this cohomology class by scaling the
$\la_r$'s in the same way.
\end{rem}
\begin{rem} \label{rem:lawe}
Although labeled projective spaces might seem to be a more general class
of toric orbifolds than orbifold projective spaces, that is not really the
case. In fact one can easily check that, up to scaling, global
coverings and/or finite quotients, the classes of labeled, orbifold and
weighted projective spaces consist of the same K\"ahler toric orbifolds.
\end{rem}

\section{Toric K\"ahler metrics} \label{sec:metrics}

Let $(M_P, \om_P, \tau_P)$ be the symplectic toric orbifold associated
to a labeled polytope $P$. In this section we describe how all
$\om_P$-compatible toric complex structures on $M_P$ (in other words,
all toric K\"ahler metrics) can be effectively parametrized by
smooth functions on $P$, according to the statements of Theorems 1 and 2.

\subsection{The ``canonical'' toric K\"ahler metric} \label{ssec:canonical}

Recall from the construction of \S\ref{ssec:toric} that
$(M_P, \om_P, \tau_P)$ comes equipped with a ``canonical'' 
$\om_P$-compatible toric complex structure $J_P$, induced from the
standard one on $\C^d$ through symplectic reduction. Following
Guillemin~\cite{Gui1},
we will now prove Theorem 1, which states that the potential $g_P$ of
$J_P$ is given by
\begin{equation}\label{eq:canpot}
g_P(x) = \frac{1}{2}\sum_{r=1}^{d} \ell_r(x)\log\ell_r(x)\,,
\end{equation}
where $\ell_r,\,r=1,\ldots,d$, are the affine functions on $(\R^n)^\ast$
defining the polytope $P$ as in~(\ref{eq:P}).

It is enough to show that the K\"ahler metric given in symplectic
$(x,\theta)$-coordinates by~(\ref{metricG}), with $G\equiv G_P = 
\Hess_x (g_P)$, corresponds to the K\"ahler metric
$\langle\cdot,\cdot\rangle_P = \om_P(\cdot,J_P\cdot)$ on $M_P$.
Because both these metrics are invariant under the $\T^n$-action and
$\breve{M}_P$ is open and dense in $M_P$, one just needs to find a suitable
slice, orthogonal to the orbits of the $\T^n$-action on $\breve{M}_P$,
and isometric to $\breve{P}$ via the moment map $\phi_P$. Here the word
``isometric'' is with respect to the metric on the slice induced by
$\langle\cdot,\cdot\rangle_P$, and the metric on $\breve{P}$ given by
$G_P$.

Such a slice arises naturally as the fixed point set of an anti-holomorphic
involution, induced from complex conjugation in $\C^d$:
\begin{equation}\label{def:inv}
\sigma:\C^d\to\C^d\,,\ \sigma(z)=\ov{z}\,,\ \Fix(\sigma)=\R^d\subset\C^d\ .
\end{equation}
The construction of \S\ref{ssec:toric} is invariant (or equivariant)
with respect to $\sigma$. In particular, the submanifold
$Z=\phi_K^{-1}(0)\subset\C^d$, with $\phi_K$ defined by~(\ref{def:phiK}),
is stable under $\sigma$. The $K$-action on $Z$ commutes with $\sigma$,
and so $\sigma$ descends to give an involution on $M_P$.

Let $Z^\si\subset\R^d$ and $M_P^\si$ denote the fixed point sets of $\si$
on $Z$ and $M_P$. Define
$$
\breve{Z}^\si = Z^\si\cap\breve{\R}^d\ \ \mbox{and}\ \ 
\breve{M}^\si_P = M_P^\si\,\cap\breve{M}_P\,,
$$
where $\breve{\R}^d=\left\{(u_1,\ldots,u_d)\in\R^d:u_r\neq 0\ \mbox{for all}
\ r=1,\ldots,d\right\}$. The following can be easily checked from the 
construction in \S\ref{ssec:toric}:
\begin{itemize}
\item[-] the quotient map $\pi:Z\to M_P$ induces a covering map
$\breve{\pi}^\si :\breve{Z}^\si\to\breve{M}_P^\si$, with group of deck
transformations given by $\{\al\in K:\al^2=1\}$;
\item[-] $\breve{\pi}^\si$ is an isometry with respect to the metric
on $\breve{Z}^\si$ induced by the Euclidean metric on $\R^d$, and the
metric on $\breve{M}_P^\si$ induced by the metric $\langle\cdot,\cdot
\rangle_P$ on $M_P$;
\item[-] the moment map $\phi_P:M_P\to P\subset (\R^n)^\ast$ induces
a covering map $\breve{\phi}_P^\si:\breve{M}_P^\si\to\breve{P}$, with
group of deck transformations given by $\{\theta\in\T^n:\theta^2=1\}$.
Moreover, $\breve{M}_P^\si$ is $\langle\cdot,\cdot\rangle_P$-orthogonal
to the orbits of the $\T^n$-action on $\breve{M}_P$.
\end{itemize}
Hence, any connected component of $\breve{M}_P^\si$ can be taken to be
the slice we were looking for. It is isometric via $\breve{\pi}^\si$
to any connected component of $\breve{Z}^\si \subset \R^d$.

Let $\breve{Z}^\si_+ = \breve{Z}^\si \cap \R^d_+$, where
$\R^d_+ = \{ (u_1,\ldots,u_d)\in\R^d : u_r >0\ \mbox{for all}\ 
r=1,\ldots,d\}$. From~(\ref{def:phiK}) we have that $\breve{Z}^\si_+$
is the subset of $\R^d_+$ defined by the quadratic equation
$$
\sum_{r=1}^{d} \left(\frac{u_r^2}{2}+\la_r\right)\iota^\ast (e_r^\ast) = 0\ .
$$
Consider the change of coordinates in $\R^d_+$ given by
$$
s_r = \frac{u_r^2}{2}\,,\ r=1,\ldots,d\,.
$$
$\breve{Z}^\si_+$ is now defined by the linear equation
$$
\sum_{r=1}^d (s_r+\la_r)\,\iota^\ast(e_r^\ast) = 0\,,
$$
and the Euclidean metric $\sum_r (du_r)^2$ becomes
\begin{equation}\label{metric-s}
\frac{1}{2}\sum_{r=1}^d \frac{(ds_r)^2}{s_r}\ .
\end{equation}
The commutative diagram~(\ref{def:phiP}) can be written here as
\begin{equation}
\begin{CD}
\breve{Z}^\si_+ @>{\breve{\phi}^\si_{\T^d}}>> (\R^d)^\ast \\
@V{\breve{\pi}^\si}VV  @AA{\beta^\ast}A \\
\breve{M}^\si_P @>{\breve{\phi}^\si_P}>>  \breve{P}\subset(\R^n)^\ast
\end{CD}
\end{equation}
and we want to determine the form of the metric~(\ref{metric-s}) on
$\breve{P}$. The map $\beta^\ast$, being dual to the surjective linear
map defined by~(\ref{def:beta}), is an injective linear map given by
\begin{equation}\label{def:beta*}
\beta^\ast (x) = \sum_{r=1}^d \langle x, m_r\mu_r\rangle\,e_r^\ast\ .
\end{equation}
The map $\breve{\phi}^\si_{\T^d}$, being the restriction of 
$\phi_{\T^d}$ defined by~(\ref{def:phiTd}), is given in the
$s$-coordinates by
\begin{equation}\label{def:phi-s}
\breve{\phi}^\si_{\T^d} (s) = \sum_{r=1}^d (s_r + \la_r)\,e_r^\ast\ .
\end{equation}
From~(\ref{def:beta*}) and~(\ref{def:phi-s}) we conclude that
$$
(s_r + \la_r) = \langle x,m_r\mu_r\rangle \Rightarrow
s_r = \langle x,m_r\mu_r\rangle - \la_r \equiv \ell_r (x)\,,\ 
\mbox{for all}\ r=1,\ldots,d\,.
$$
Hence the metric~(\ref{metric-s}) can be written in the $x$-coordinates
of the polytope $P$ as
$$
\frac{1}{2}\sum_{r=1}^d \frac{(ds_r)^2}{s_r} =
\frac{1}{2}\sum_{r=1}^d \frac{(d\ell_r)^2}{\ell_r} =
\sum_{i,j=1}^n \frac{\p^2 g_P (x)}{\p x_i \p x_j}\,dx_i dx_j\,,
$$
where $g_P$ is given by~(\ref{eq:canpot}) and the last equality
is a trivial exercise.

This completes the proof of Theorem~\ref{thm1}.

\subsection{General toric K\"ahler metrics} \label{ssec:general}

We will now prove Theorem~\ref{thm2}, which states that on a symplectic toric
orbifold $(M_P,\om_P,\tau_P)$, associated to a labeled polytope $P$,
compatible toric complex structures $J$ are in one to one correspondence
with potentials $g\in C^\infty({\breve{P}})$ of the form
\begin{equation} \label{def:gen-g}
g = g_P + h\,,
\end{equation}
where $g_P$ is given by~(\ref{eq:canpot}), $h$ is smooth on the whole $P$,
and the matrix $G = \Hess (g)$ is positive definite on $\breve{P}$ and has
determinant of the form
\begin{equation} \label{eq:det-G}
\Det (G) = \left[\de \prod_{r=1}^d \ell_r \right]^{-1}\,,
\end{equation}
with $\de$ being a smooth and strictly positive function on the whole $P$.

The proof of this theorem for symplectic toric orbifolds given in the
Appendix of~\cite{Abr2}, generalizes with very minor modifications
to our orbifold context.

We first prove that any potential $g\in C^\infty (\breve{P})$ of the
form~(\ref{def:gen-g}) and satisfying~(\ref{eq:det-G}), defines
through~(\ref{matrixJ}) a compatible toric complex structure $J$
on $(M_P,\om_P,\tau_P)$. It is clear that $J$ is well defined on
$\breve{M}_P \cong \breve{P}\times\T^n$. To see that it extends to the
whole $M_P$ one has to check that the singular behaviour of $J$ near
the boundary of $P$ is the same as the singular behaviour of $J_P$,
which we know extends to the whole $M_P$.

This singular behaviour is best described in terms of the Hessians
$G_P=\Hess(g_P)$ and $G=\Hess(g)$. Explicit calculations show that
although $G_P$ is singular on the boundary of the polytope $P$, 
$G_P^{-1}$ is smooth on the whole $P$ and its determinant has the form
$$
\Det (G_P^{-1}) = \de_P \prod_{r=1}^d \ell_r\,,
$$
where $\de_P$ is a smooth and strictly positive function on the
whole $P$. This formula captures the relevant singular behaviour and
has the following geometric interpretation. Given $x\in P$,
let $\Ff(x)$ be the set of facets of $P$ that contain $x$, i.e.
$$
\Ff(x) = \left\{ r\in\{1,\ldots,d\}:\ell_r(x) = 0 \right\}\ .
$$
The kernel of $G_P^{-1}(x)$ is precisely the linear span of the
normals $\mu_r\in\R^n$ for $r\in\Ff(x)$. Due to Remark~\ref{rem:isorb},
this kernel is also the Lie algebra of the isotropy group
$\Ga_p\subset\T^n$ of any $p\in M_P$ such that $\phi_P(p)=x$.
Conditions~(\ref{def:gen-g}) and~(\ref{eq:det-G}) guarantee that
$G^{-1}$ has these same degeneracy properties, and that is enough for
the corresponding $J$ to extend to a compatible toric complex structure
well defined on the whole $M_P$.

We now prove that any compatible toric complex structure $J$ on
$(M_P,\om_P,\tau_P)$ corresponds, in suitable symplectic coordinates
on $\breve{M}_P \cong \breve{P}\times\T^n$, to a potential
$g\in C^\infty(\breve{P})$ of the form~(\ref{def:gen-g}). Because
$J$ is apriori defined on the whole $M_P$, the matrix $G=\Hess (g)$
will automatically satisfy~(\ref{eq:det-G}). The idea of the proof
is to translate to symplectic coordinates some well known facts
from K\"ahler geometry.

It follows from Theorem~9.4 in~\cite{LeTo} that there is an equivariant
biholomorphism
$$
\varphi_J:(M_P, J_P, \tau_P) \to (M_P, J, \tau_P)\,,
$$
with $\varphi_J$ acting as the identity in cohomology. This means that
$(M_P,\om_P,J)$ is equivariantly K\"ahler isomorphic to 
$(M_P,\om_J,J_P)$, where $\om_J = (\varphi_J)^\ast (\om_P)$ and
$[\om_J] = [\om_P] \in H^2 (M_P)$. It follows from~\cite{Bail} that the
$\p\ov{\p}$-lemma is valid on K\"ahler orbifolds, and hence there exists
a $\T^n$-invariant smooth function $f_J \in C^\infty (M_P)$ such that
$$
\om_J = \om_P + 2i\p\ov{\p}f_J\,,
$$
where the $\p$- and $\ov{\p}$-operators are defined with respect to the
complex structure $J_P$.

In the symplectic $(x,\theta)$-coordinates on $\breve{M}_P\cong
\breve{P}\times\T^n$, obtained via the ``canonical'' moment map
$\phi_P : M_P \to P\subset (\R^n)^\ast$, we then have a function
$f_J \equiv f_J(x)$, smooth on the whole polytope $P$, and such that
$$
\om_J = dx\wedge d\theta + 2i\p\ov{\p}f_J\ .
$$
The rest of the proof consists of the following three steps:
\begin{itemize}
\item[(i)] write down on $P$ the change of coordinates
$\tvphi_J:P\to P$ that corresponds to the equivariant biholomorphism
$\varphi_J:M_P\to M_P$ and transforms the symplectic $(x,\theta)$-coordinates
for $\om_P$ into symplectic $(\tx=\tvphi_J(x),\theta)$-coordinates
for $\om_J$. These $(\tx,\theta)$-coordinates are the suitable symplectic
coordinates we were looking for;
\item[(ii)] find the potential $g=g(\tx)$ for the transformed
$J=(\tvphi_J)_\ast (J_P)$ in these $(\tx,\theta)$-coordinates;
\item[(iii)] check that the function $h:\breve{P}\to\R$ given by
$h(\tx)=g(\tx)-g_P(\tx)$, with $g_P(\tx) = \frac{1}{2}\sum_r
\ell_r(\tx)\log\ell_r(\tx)$, is actually defined and smooth on the 
whole $P$. Here, as always, $\ell_r(\tx)\equiv \langle\tx,m_r\mu_r\rangle
-\la_r\,,\ r=1,\ldots,d$, are the defining functions of the polytope $P$.
\end{itemize}

All these steps can be done in a completely explicit way. We refer the
reader to the Appendix in~\cite{Abr2} for details. The change of 
coordinates in step (i) is given in vector form by
$$
\tx=\tvphi_J(x) = x + G_P^{-1}\cdot\frac{\p f_J}{\p x}\,,
$$
where $\p f_J/\p x = (\p f_J/\p x_1,\ldots,\p f_J/\p x_n)^t\equiv$
column vector. $\tvphi_J$ is a diffeomorphism of the whole $P$ and,
due to the degeneracy behaviour of the matrix $G_P^{-1}$ on the 
boundary of $P$, preserves each of its faces (i.e. each vertex, edge,
\dots,facet and interior $\breve{P}$). In step (ii) one finds that
$$
g(\tx) = \left\langle \tx - \tvphi_J^{-1}(\tx), (\frac{\p g_P}{\p x}\circ
\tvphi_J^{-1})(\tx)\right\rangle + (g_P\circ\tvphi_J^{-1})(\tx)
- (f_J\circ\tvphi_J^{-1})(\tx)\ .
$$
Since $(f_J\circ\tvphi_J^{-1}) \in C^\infty (P)$ and $\tvphi_J$ is a
smooth diffeomorphism of the whole $P$, step (iii) reduces to
checking that
$$
\left\langle \tvphi_J(x) -x , \frac{\p g_P}{\p x}(x)\right\rangle
+ g_P(x) - g_P(\tvphi_J(x)) \in C^\infty (P)\ .
$$
Simple explicit computations show that this is true provided
$$
\frac{\ell_r(x)}{\ell_r(\tvphi_J(x))} \in C^\infty (P)\  
\mbox{for all}\ \ r=1,\ldots,d\,,
$$
and this follows from the fact that $\tvphi_J$ preserves the
combinatorial structure of $P$.

The proof of Theorem~\ref{thm2} is completed.

\section{Extremal metrics} \label{sec:extremal}

In this section, after some preliminaries on extremal K\"ahler metrics,
we use the framework of Sections~\ref{sec:symp} and~\ref{sec:metrics}
to prove Theorem~\ref{thm3}, i.e. we give a simple description of a toric
extremal K\"ahler metric on any labeled projective space. Due to
Proposition~\ref{prop:orbproj} and Remark~\ref{rem:cover}, this gives
rise in particular to a toric extremal K\"ahler metric on any weighted
projective space, as stated in Corollary~\ref{cor1}. Theorem~\ref{thm4}
is proved in the last subsection.

\subsection{Preliminaries on extremal metrics} \label{ssec:pre-ext}

In~\cite{Cal1} and~\cite{Cal2}, Calabi introduced the notion of {\bf extremal}
K\"ahler metrics. These are defined, for a fixed closed complex manifold
$(M,J_0)$, as critical points of the square of the $L^2$-norm of the scalar
curvature, considered as a functional on the space of all symplectic K\"ahler
forms $\om$ in a fixed K\"ahler class $\Om\in H^2(M,\R)$. The extremal
Euler-Lagrange equation is equivalent to the gradient of the scalar curvature
being an holomorphic vector field (see~\cite{Cal1}), and so these metrics
generalize constant scalar curvature K\"ahler metrics. Calabi illustrated
this in~\cite{Cal1} by constructing families of extremal K\"ahler metrics
of non-constant scalar curvature. Moreover, Calabi
showed in~\cite{Cal2} that extremal K\"ahler metrics are always invariant under
a maximal compact subgroup of the group of holomorphic transformations of
$(M,J_0)$. Hence, on a complex toric manifold or orbifold, 
extremal K\"ahler metrics are automatically toric K\"ahler metrics, 
and one should be able to write them down using the framework of
Section~\ref{sec:metrics}. This was carried out in~\cite{Abr1} for
Calabi's simplest family, having $\cp\#\,\ocp$ as underlying toric manifold.

We now recall from~\cite{Abr1} some relevant differential-geometric formulas
in symplectic $(x,\theta)$-coordinates. A K\"ahler metric of the
form~(\ref{metricG}) has scalar curvature $S$ given by\footnote[1]{The
normalization for the value of the scalar curvature we are using here is 
the same as in~\cite{Bess}. It differs from the one used in~\cite{Abr1,Abr2} 
by a factor of $1/2$.}
\begin{equation} \label{scalarsymp1}
S = - \sum_{j,k} \frac{\p}{\p x_j}
\left( g^{jk}\, \frac{\p \log \Det(G)}{\p x_k} \right)\,,
\end{equation}
which after some algebraic manipulations becomes the more compact
\begin{equation} \label{scalarsymp2}
S = - \sum_{j,k} \frac{\p^2 g^{jk}}{\p x_j \p x_k}\,, 
\end{equation}
where the $g^{jk},\ 1\leq j,k\leq n$, are the entries of the inverse 
of the matrix $G = \Hess_x (g)$, $g\equiv$ potential. 
The Euler-Lagrange 
equation defining an extremal K\"ahler metric can be shown to be equivalent to
\begin{equation} \label{extremalsymp}
\frac{\p S}{\p x_j} \equiv\ \mbox{constant},\ j=1,\ldots,n,
\end{equation}
i.e. {\it the metric is extremal if and only if its scalar curvature $S$ is an
affine function of $x$}. One can express~(\ref{extremalsymp}) in more invariant
terms, giving a symplectic analogue of the complex extremal condition saying
that the gradient of the scalar curvature is an holomorphic vector field.
\begin{prop} \label{prop:extsymp}
Let $(M_P,\om_P,\tau_P)$ be a compact symplectic toric orbifold with
moment map $\phi_P:M_P\to P\subset (\R^n)^\ast$. A toric compatible
complex structure $J$ gives rise to an extremal K\"ahler metric $\langle\cdot,
\cdot\rangle = \om_P(\cdot,J\cdot)$ if and only if its scalar curvature
$S$ is a constant plus a linear combination of the components of the
moment map $\phi_P$.

In other words, the metric is extremal if and only if there exists
$\xi\in\R^n\equiv$ Lie algebra of $\T^n$, such that
$$
dS = d\langle\xi,\phi_P\rangle\ .
$$
\end{prop}

\subsection{Extremal orbifold metrics on $S^2$} \label{ssec:ext-S2}

Here we prove Theorem~\ref{thm3} when $n=1$. This very simple case
is already interesting and motivates well the formula for the
potential $g$ in the general case.

Consider the one dimensional labeled polytope defined by
$$
\ell_1(x) = m_1 (1+x)\ \ \mbox{and}\ \ \ell_2(x) = m_2 (1-x)\,,
\ \mbox{with}\ m_1,m_2\in\N\,.
$$
The corresponding labeled projective space $\lp^2_\bm$ is homeomorphic
to the $2$-sphere $S^2$, and the orbifold structure at each pole can be
geometrically interpreted as a conical singularity with angle
$2\pi/m_r\,,\ r=1,2$. We look for an extremal metric generated by a
potential $g\in C^\infty(-1,1)$ of the form
$$
g(x) = \frac{1}{2}\left( m_1(1+x)\log(m_1(1+x)) + m_2(1-x)\log(m_2(1-x))
+ h(x) \right)\,,
$$
with $h\in C^\infty[-1,1]$. Formula~(\ref{scalarsymp2}) for the scalar
curvature becomes
$$
S(x) = - \left(\frac{1}{g''(x)}\right)'' =
\left(-(1-x^2)\ov{h}(x)\right)''
$$
where
$$
\ov{h}(x) = \frac{2}{m_1(1-x)+m_2(1+x) + (1-x^2)h''(x)} \in C^\infty[-1,1]\,.
$$
Equation~(\ref{extremalsymp}) says that the metric is extremal if and
only if $S$ is a first degree polynomial, hence if and only if $\ov{h}$
is a first degree polynomial. Since $\ov{h}(-1)=1/m_1$ and $\ov{h}(1)=1/m_2$
we must have
$$
\ov{h}(x) = \frac{1}{2m_2}(1+x) + \frac{1}{2m_1}(1-x) =
\frac{\ell_1(x) + \ell_2(x)}{2m_1 m_2}\,.
$$
Solving for $h''(x)$ and integrating one gets
$$
h(x) = - (m_1(1+x) + m_2(1-x))\log (m_1(1+x) + m_2(1-x))\,,
$$
i.e.
$$
h = -\ell_\Si\log\ell_\Si\ \ \mbox{with}\ \ \ell_\Si = \ell_1 + \ell_2\,.
$$
Note that, because $\ell_\Si$ is strictly positive on $[-1,1]$, $h$ is
defined and smooth on $[-1,1]$. Moreover, 
$$
G^{-1} = \frac{1}{g''} = \frac{\ell_1\,\ell_2\,\ell_\Si}{2\,m_1^2\,m_2^2}
$$
is strictly positive on $(-1,1)$ and has the degeneracy behaviour at
the boundary points $-1$ and $1$ required by~(\ref{eq:det-G}). 

Hence the potential
\begin{equation}\label{potext1}
g = \frac{1}{2}\left(\ell_1\log\ell_1 +\ell_2\log\ell_2 - \ell_\Si\log
\ell_\Si\right)
\end{equation}
defines a toric extremal K\"ahler metric on $\lp^2_\bm$. Its scalar
curvature is given by
$$
S(x) = \frac{(m_1 + m_2) + 3x(m_1 - m_2)}{m_1 m_2}\,.
$$
As a function on $\lp^2_\bm$ it can be written as
$$
S = \left(\frac{1}{m_1} + \frac{1}{m_2}\right) + 3
\left(\frac{1}{m_2}-\frac{1}{m_1}\right)\,\phi_\bm\,,
$$
where $\phi_\bm:\lp^2_\bm \to [-1,1]\subset\R^\ast$ is the moment map.
Hence
$$
dS = d\langle\xi_\bm,\phi_\bm\rangle\ \ \mbox{for}\ \ 
\xi_\bm= 3\left(\frac{1}{m_2}-\frac{1}{m_1}\right)
\in\R\cong\mbox{Lie algebra of}\ \T^1\,.
$$

\subsection{Extremal metrics on $\lp^n_\bm$} \label{ssec:ext-SP}

We now consider a general labeled simplex $P^n_\bm\subset(\R^n)^\ast$
defined by
$$
\ell_r(x) = m_r (1+x)\,,\ r=1,\ldots,n\,,\ \ell_{n+1}(x) =
m_{n+1}(1-\psi)\,,\ \psi = \sum_{j=1}^n x_j\,,
$$
with $m_r\in\N$, for all $r=1,\ldots,n+1$. The corresponding labeled
projective space $\lp^n_\bm$ is homeomorphic to $\cp^n$ (see 
Remark~\ref{rem:spcp}). Under this homeomorphism the pre-image of
the $r$-th facet
$$
F_r = \left\{ x\in P^n_\bm: \ell_r(x) = 0 \right\}
$$
by the moment map $\phi_\bm : \lp^n_\bm \to P^n_\bm$ corresponds to
$$
N_r = \left\{ [z_1,\ldots,z_{n+1}]\in\cp^n : z_r = 0\right\}
\cong \cp^{n-1}\,.
$$
The orbifold structure of $\lp^n_\bm$ can be geometrically interpreted
on $\cp^n$ as a conical singularity with angle $2\pi/m_r$ around each
$N_r\cong\cp^{n-1}$, for $r=1,\ldots,n+1$.

Motivated by the form of the potential~(\ref{potext1}) for the toric extremal
K\"ahler metric on $\lp^2_\bm$, we consider now the potential
$g\in C^\infty(\breve{P}^n_\bm)$ given by
\begin{equation}\label{potextn}
g = \frac{1}{2}\left(\sum_{r=1}^{n+1} \ell_r\log\ell_r
- \ell_\Si\log\ell_\Si\right)\ \ \mbox{with}\ \ 
\ell_\Si = \sum_{r=1}^{n+1} \ell_r\ .
\end{equation}
Note that since $\ell_\Si$ is strictly positive on $P^n_\bm$, this
potential $g$ is of the general form~(\ref{def:gen-g}). 

The entries of the matrix $G = \Hess (g)$ are easily computed 
from~(\ref{potextn}):
\begin{equation} \label{eq:m-gjk}
g_{jk} = \frac{\p^2 g}{\p x_j \p x_k} = \frac{1}{2}\left(
\delta_{jk}\frac{m_j^2}{\ell_j} + \frac{m_{n+1}^2}{\ell_{n+1}} -
\frac{(m_j - m_{n+1})(m_k - m_{n+1})}{\ell_\Si}\right)\,,
\end{equation}
where $\delta_{jk}$ is equal to $1$ if $j=k$ and equal to $0$ otherwise.
The proof of the following lemma is left as an exercise to the reader.
\begin{lem} \label{lem:detinv}
The matrix $G = \Hess (g) = (g_{jk})_{j.k=1}^{n,n}$ is positive definite
on $\breve{P}^n_\bm$ with determinant given by
\begin{equation} \label{eq:det-spnm}
\Det(G) = \left[ \left(\prod_{r=1}^{n+1}\ell_r\right)
\frac{2^n \ell_\Si}{(n+1)^2 \prod_{r=1}^{n+1} m_r^2}\right]^{-1}\ .
\end{equation}
The entries of the matrix $G^{-1} = (g^{jk})_{j,k=1}^{n,n}$ are given by
\begin{equation} \label{eq:inv-spnm}
g^{jk} = 2 \left(\delta_{jk}\frac{\ell_j}{m_j^2} - \frac{m_j + m_k}{n+1}
\frac{\ell_j\,\ell_k}{m_j^2\,m_k^2} + \frac{1}{(n+1)^2}
\frac{\ell_j\,\ell_k}{m_j\,m_k}\left(\sum_{r=1}^{n+1}\frac{\ell_r}{m_r^2}
\right)\right)\,.
\end{equation}
\end{lem}

It follows that the potential $g$ defined by~(\ref{potextn}) satisfies
the conditions of Theorem~\ref{thm2}, and hence defines a toric K\"ahler metric
on $\lp^n_\bm$. Moreover, since each $g^{jk}$ is a third degree polynomial,
it is clear from~(\ref{scalarsymp2}) that the scalar curvature $S$ is
a first degree polynomial. By~(\ref{extremalsymp}) this means that the
metric defined by $g$ is indeed extremal, thus finishing the proof of
Theorem~\ref{thm3}.

More explicitly, we have that the scalar curvature is given by
$$
S(x) = \frac{2n}{n+1}\left(\sum_{r=1}^{n+1} \frac{1}{m_r}\right)
+ \frac{2(n+2)}{n+1} \sum_{j=1}^n \left(\frac{1}{m_{n+1}}-\frac{1}{m_j}
\right)\,x_j\ .
$$
As a function on $\lp^n_\bm$ it can be written as
$$
S =  \frac{2n}{n+1}\left(\sum_{r=1}^{n+1} \frac{1}{m_r}\right) +
\langle \xi_\bm,\phi_\bm\rangle\,,
$$
where $\phi_\bm$ is the moment map and
$$
\xi_\bm = \frac{2(n+2)}{n+1} \left(\frac{1}{m_{n+1}}-\frac{1}{m_1},
\ldots,\frac{1}{m_{n+1}}-\frac{1}{m_n}\right)\in \R^n \cong
\mbox{Lie algebra of}\ \T^n\ .
$$

\subsection{Conical extremal metrics on $\lp^n_\bo$}\label{ssec:conical}

The purpose of this subsection is to prove Theorem~\ref{thm4}, i.e. we
will describe natural ``conical'' compactifications of extremal K\"ahler
metrics defined by potentials $g$ of the form~(\ref{potextn}), for
any positive real vector of labels $\bm\in\R^{n+1}_+$.

The symplectic toric orbifold where this compactification takes place
is obtained by forgetting the labels. Hence we consider the standard
smooth symplectic toric manifold $(\lp^n_\bo, \om_\bo, \tau_\bo)$
associated to the simplex $P^n_\bo\subset(\R^n)^\ast$, and denote by
$\phi_\bo : \lp^n_\bo\to P^n_\bo$ the corresponding moment map. Note
that $(\lp^n_\bo, \om_1, \tau_1)$ is equivariantly symplectomorphic to
$\cp^n$ with a suitably normalized Fubini-Study symplectic form and
standard torus action.

For any $m\in\R^{n+1}_+$, the potential $g\in C^\infty(\breve{P}^n_\bo)$
given by~(\ref{potextn}) defines an extremal K\"ahler metric
$\langle\cdot,\cdot\rangle_\bm$ on $\breve{\lp}^n_\bo = \phi_\bo^{-1}
(\breve{P}_\bo^n) \cong \breve{P}^n_\bo \times \T^n$ given by~(\ref{metricG}).
Consider the pre-image $N_r \equiv \phi_\bo^{-1} (F_r)$ of each facet
$F_r\subset P_\bo^n\,,\ r=1,\ldots,n+1$. Each $N_r$ is a real codimension
$2$ symplectic toric submanifold of $\lp^n_\bo$, symplectomorphic to a
suitably normalized $\lp^{n-1}_\bo \cong \cp^{n-1}$. The restriction
$\phi_\bo |_{N_r} : N_r \to F_r$ is a corresponding moment map. We want
to show that $\langle\cdot,\cdot\rangle_\bm$ extends to an extremal metric
on the whole $\lp^n_\bo$ with conical singularities of angles $2\pi/m_r$
around each $N_r$.

The potential $g$, although only smooth on the interior $\breve{P}^n_\bo$,
is a continuous function on the whole polytope $P^n_\bo$. Denote by
$g_r\in C^\infty(\breve{F}_r)\cap C^0 (F_r)$ the restriction of $g$ to
$F_r$ (here $\breve{F}_r$ denotes the relative interior of $F_r$).
Using the explicit form of the matrix $G=\Hess(g)$ given by~(\ref{eq:m-gjk}),
one can easily check that the extremal metric $\langle\cdot,\cdot\rangle_\bm$,
defined on $\breve{\lp}^n_\bo$, extends to a well defined smooth extremal 
metric on $\breve{N}_r \equiv \phi_\bo^{-1}(\breve{F}_r)$ whose potential
is exactly given by $g_r$. Note that the hyperplane in $(\R^n)^\ast$ that
contains $F_r$ has an induced affine structure, and so it makes sense to
consider $G_r = \Hess (g_r)$.

Because of the equivariant version of Darboux's theorem, we can understand
what happens in the normal directions to each point $p\in\breve{N}_r$ by
analysing a neighborhood of zero in $\R^2$. In $(r,\theta)$-polar
coordinates the standard symplectic form is $rdr\wedge d\theta$, and the
moment map for the standard circle action is given by $x=r^2/2$. The moment
polytope is $[0,+\infty)$ defined by the single affine function
$\ell(x)=x$. The standard smooth K\"ahler metric is defined by the
potential $g_1 = \frac{1}{2} x\log x$, hence given by
$$
\langle \cdot, \cdot \rangle_1 = G_1''\,dx^2 + \frac{1}{g_1''}\, d\theta^2 = 
\frac{1}{2x}\,dx^2 + 2x\, d\theta^2\,,
$$
while the ``orbifold'' one is defined for any $m\in\R_+$ by the potential
$g_m = \frac{1}{2} mx \log (mx)$, and hence given by
$$
\langle \cdot, \cdot \rangle_m = g_m''\,dx^2 + \frac{1}{g_m''}\, d\theta^2 = 
\frac{m}{2x}\,dx^2 + \frac{2x}{m}\, d\theta^2\,.
$$
In $(r,\theta)$-polar coordinates we get
$$
\langle \cdot, \cdot \rangle_1 = dr^2 + r^2 d\theta^2 \equiv
\ \mbox{standard smooth flat metric}\,,
$$
while
$$
\langle \cdot, \cdot \rangle_m = m\left( dr^2 + \left(\frac{r}{m}\right)^2
d\theta^2 \right)
$$
which is the polar form of a metric with a conical singularity of angle
$2\pi/m$ around the origin.

Hence we have an extension of each extremal K\"ahler metric
$\langle\cdot,\cdot\rangle_\bm\,,\ 
\bm\in\R^{n+1}_+$, from $\breve{\lp}^n_\bo$
to $\breve{\lp}^n_\bo \cup \left( \cup_{r=1}^{n+1} \breve{N}_r\right)$,
having normal conical singularities around each $\breve{N}_r$. The same
argument can be used to show that the metric on each $\breve{N}_r$ extends
to the moment map pre-images of the relative interior of each facet of
$F_r$ (an $(n-2)$-dimensional simplex and face of $P^n_\bo$). One can
continue this process until the metric is extended to the whole
$\lp^n_\bo$. For example, at the last step one extends the metric to the
fixed points of the $\T^n$-action, corresponding to the vertices of
$P^n_\bo$. There the metric looks like the product of $n$ 
cones of dimension two and angles $2\pi/m_{r_i}$,              
where $m_{r_1},\ldots,m_{r_n}$ are the
positive real labels of the $n$ facets of $P^n_\bo$ that meet at the 
relevant vertex.

\section{A family of self-dual Einstein metrics} \label{sec:sde}

Recall from the introduction that the extremal K\"ahler metrics given
by Theorem~\ref{thm4} are actually Bochner-K\"ahler (see~\cite{Brya}).
In dimension four ($n=2$) ``Bochner-K\"ahler'' is the same as
``self-dual K\"ahler''. It follows from the work of Derdzinski~\cite{Derd}
and Apostolov-Gauduchon~\cite{ApGa} that, whenever its scalar curvature
$S$ is nonzero, a self-dual K\"ahler metric is conformally Einstein, with
conformal factor given by $S^{-2}$. In this section we explore this relation
for a particular one-parameter family of metrics arising from 
Theorem~\ref{thm4}.

Consider the smooth symplectic toric manifold $(\lp^2_\bo\cong\cp^2,
\om_\bo,\tau_\bo)$ associated to the simplex $P^2_\bo\subset (\R^2)^\ast$.
For any $\bm=(1,1,m)\,,\ m\in\R_+$, let $\langle\cdot,\cdot\rangle_m$ be
the extremal K\"ahler metric defined by the potential
\begin{equation} \label{eq:potgm}
g_m(x) = \frac{1}{2}\left(\sum_{r=1}^3 \ell_r(x) \log\ell_r(x) -
\ell_\Si (x) \log \ell_\Si (x) \right)\,,
\end{equation}
where $\ell_1(x) = 1+x_1\,,\ \ell_2(x) = 1+x_2\,,\ 
\ell_3(x) = m(1-\psi)$ and $\ell_\Si(x) = 2+m - (m-1)\psi$.
Here and in the rest of this section $\psi=x_1+x_2$. Note that 
in~(\ref{eq:potgm}) the two terms with $\ell_1$ and $\ell_2$
correspond to the standard flat metric on $\R^4$, while the terms
with $\ell_3$ and $\ell_\Si$ only depend on the ``radial'' coordinate
$\psi$. This means that the metric $\langle\cdot,\cdot\rangle_m$
defined by the potential $g_m$ is $U(2)$-invariant (see~\cite{Abr3} for
a general discussion of this type of metrics).

The scalar curvature $S_m$ of $\langle\cdot,\cdot\rangle_m$ is given by
\begin{equation}\label{eq:scm}
S_m (x) = \frac{4}{3m}(2m+1 + 2(1-m)\psi)\,,
\end{equation}
which is strictly positive on $\lp^2_\bo$ if $m>1/2$. Hence, for any
$1/2 < m < +\infty$, the metric $\langle\cdot,\cdot\rangle_{m}^\ast
\equiv S_m^{-2}\langle\cdot,\cdot\rangle_{m}$ is a self-dual Einstein
metric on $\lp^2_\bo\cong \cp^2$ with a normal conical singularity
of angle $2\pi/m$ around a $\lp^1_\bo \cong \cp^1$.

In this simple case it is not hard to check explicitly that
$\langle\cdot,\cdot\rangle_{m}^\ast$ is Einstein and compute its
scalar curvature. A result of Derdzinski~\cite{Derd} (see also~\cite{Bess})
states that for any $4$-dimensional extremal K\"ahler metric 
$\langle\cdot,\cdot\rangle$ with non constant scalar curvature $S$, the
metric $\langle\cdot,\cdot\rangle^\ast\equiv S^{-2}\langle\cdot,\cdot\rangle$ 
is Einstein if and only if
\begin{equation} \label{eq:extEin}
S^3 - 6S\Delta S - 12 |dS|^2 = \ \mbox{constant}\,.
\end{equation}
Moreover, a standard formula for the scalar curvatures of conformally
related metrics (see e.g.~\cite{Bess}) states that the scalar curvature
$S^\ast$ of $\langle\cdot,\cdot\rangle^\ast$ is given by
\begin{equation} \label{eq:confsc}
S^\ast = S^3 (6\Delta(S^{-1}) + 1)\,.
\end{equation}
In both these formulas $\Delta$ is the Laplacian with respect to the metric 
$\langle\cdot,\cdot\rangle$.

For any toric K\"ahler metric defined by a potential $g\in C^\infty
(\breve{P})$, the Laplacian $\Delta$ of a function $f\in C^\infty(P)$
(i.e. a smooth $\T^n$-invariant function on $M_P$) is given by
$$
\Delta f = -(\Det G)\sum_{j,k=1}^n g^{jk} \frac{\p}{\p x_j}
\left(\frac{1}{\Det G}\frac{\p f}{\p x_k}\right)\,,
$$
where $G = \Hess(g) = (g_{jk})$ and $g^{jk}$ are the entries of $G^{-1}$.
This formula, together with the simple form of $g_m$ and $S_m$, makes
the calculations involved in~(\ref{eq:extEin}) and~(\ref{eq:confsc})
easy enough.

For example, one computes that the scalar curvature $S_m^\ast$ of
$\langle\cdot,\cdot\rangle_{m}^\ast$ is given by
$$
S_m^\ast = \left(\frac{4}{m}\right)^3 (2m-1)\,.
$$
One sees that the self-dual Einstein metric 
$\langle\cdot,\cdot\rangle_{m}^\ast$ has positive scalar curvature when
$1/2 <m<+\infty$, but is actually Ricci-flat when $m=1/2$ or $m=+\infty$
provided we can make sense of it.

When $m=1/2$ the extremal scalar curvature $S_{1/2}$ is given by
$$
S_{1/2} (x) = \frac{8}{3} (2+\psi)\,,
$$
and hence vanishes at the unique point of $\lp^2_\bo\cong\cp^2$ corresponding
to the vertex $(-1,-1)\in P^2_\bo$. The complement of this point in
$\cp^2$ is just the normal bundle of the ``opposite'' $\cp^1$ (corresponding
to the facet $F_3\subset P^2_\bo$), i.e. a line bundle with first Chern
class $c_1=1$. The label $m=1/2$ means that the normal conical singularity
can be resolved by passing to a $\Z_2$-quotient, i.e. to the line bundle
with $c_1=2$ given by $T\,\cp^1$. This means that the self-dual 
Ricci-flat Einstein metric $\langle\cdot,\cdot\rangle_{1/2}^\ast=
S_{1/2}^{-2}\langle\cdot,\cdot\rangle_{1/2}$ is smooth and complete when
considered on $T\,\cp^1$. Being $U(2)$-invariant, it must coincide with
the well-known Eguchi-Hanson metric~\cite{EgHa}.

When $m\to\infty$ the matrix $G_m = \Hess (g_m)$ converges to the matrix
$$
G_\infty (x) = \frac{1}{2}
\begin{bmatrix}
\frac{1}{1+x_1} + \frac{4-\psi}{(1-\psi)^2} & \frac{4-\psi}{(1-\psi)^2} \\
\ &\ \\
\frac{4-\psi}{(1-\psi)^2} & \frac{1}{1+x_2} + \frac{4-\psi}{(1-\psi)^2}
\end{bmatrix}
\,.
$$
One easily checks that $G_\infty = \Hess (g_\infty)$ where
$$
g_\infty (x) = \frac{1}{2}\left( \sum_{r=1}^3 \ell_r(x)\log \ell_r(x)
- 3\log (1-\psi)\right)\,.
$$
The metric $\langle\cdot,\cdot\rangle_\infty$ defined by this potential
does not extend to the whole $\lp^2_\bo\cong\cp^2$. However it is 
a well-defined smooth complete extremal K\"ahler metric of finite volume
on $B=\lp^2_\bo\setminus\lp^1_\bo$, where the sphere $\lp^1_\bo\cong\cp^1$
corresponds to the facet $F_3\subset P^2_\bo$. In the normal directions
to this sphere at infinity, the extremal metric 
$\langle\cdot,\cdot\rangle_\infty$ looks like a complete hyperbolic cusp
(this can be seen by considering for example $m_1=1$ and $m_2\to +\infty$
for the orbifold metrics on $S^2$ discussed in~\S\ref{ssec:ext-S2}).
Note that $B$ is symplectomorphic to an open ball in $\R^4$ and,
with respect to the complex structure $J_\infty$ defined by $g_\infty$,
biholomorphic to $\C^2$.

The scalar curvature of $\langle\cdot,\cdot\rangle_\infty$ is given by
$$
S_\infty (x) = \frac{8}{3} (1-\psi)\,,
$$
which vanishes exactly at the sphere at infinity. Hence, the metric
$\langle\cdot,\cdot\rangle_\infty^\ast \equiv S_\infty^{-2}\,
\langle\cdot,\cdot\rangle_\infty$ is a smooth complete self-dual
Ricci-flat Einstein metric on $B$, obviously with infinite volume.
Being $U(2)$-invariant and $B$ being diffeomorphic to $\R^4$, it
must coincide with the well-known Taub-NUT metric~\cite{EGHa}.

As promised in the introduction, we get in this way a one parameter
family of $U(2)$-invariant self-dual Einstein metrics
$\langle\cdot,\cdot\rangle_m^\ast$, $1/2\leq m\leq +\infty$,
having positive scalar curvature when $1/2 < m < +\infty$ and
connecting the Ricci-flat Eguchi-Hanson metric on $T\,\cp^1$ 
$(m=1/2)$ with the Ricci-flat Taub-NUT metric on $\R^4\cong\C^2$
$(m=+\infty)$. Note that one of the metrics in between is the 
K\"ahler-Einstein Fubini-Study metric on $\cp^2$ $(m=1)$.

In~\cite{GaLa} Galicki and Lawson use quaternionic reduction to
produce self-dual Einstein metrics on certain weighted projective spaces.
These include $\cp^2_{(p+q,p+q,2p)}$, which up to covering/quotient
correspond in the above family to $m=(p+q)/2p$. Galicki-Lawson assume that
$p,q\in\N$, $q\leq p$ and $(p,q)=1$. They point out that when $q/p\to 1$
their metrics converge to Fubini-Study on $\cp^2$, while when $q/p\to 0$
they converge to Eguchi-Hanson on $T\,\cp^1$. This is consistent with the
$m=1$ and $m=1/2$ cases in our family. In fact, it follows from the
classification results of~\cite{ApGa} that the Galicki-Lawson metrics,
whenever defined, are the same as the ones constructed here for the
corresponding value of the parameter $m$.


\subsection*{Acknowledgments}

I would like to thank the support and hospitality of The Fields Institute
for Research in Mathematical Sciences, where this work was carried out,
in particular the organizers of the program in Symplectic Topology, Geometry 
and Gauge Theory (January-June, 2001): Lisa Jeffrey, Boris Khesin and
Eckhard Meinrenken.

I would also like to thank Vestislav Apostolov, Andrew Dancer, Paul Norbury 
and Susan Tolman for helpful conversations regarding this paper.


\end{document}